\documentclass[twocolumn]{autart}    

\usepackage{graphicx}          

\usepackage[figuresright]{rotating}
\newtheorem{theorem}{Theorem}
\newtheorem{remark}{Remark}
\newtheorem{lemma}{Lemma}
\newtheorem{proposition}{Proposition}

\newtheorem{assumption}{Assumption}
\usepackage{bbm}
\usepackage{graphics}
\usepackage{epsfig}
\usepackage{mathptmx}
\usepackage{times}
\usepackage{amsmath}
\usepackage{amssymb}
\usepackage{float}
\usepackage{xcolor}

\usepackage{cite}
\usepackage{subfig}
\usepackage{psfrag}
\usepackage{url}
\usepackage{stfloats}

\usepackage{algorithm}
\usepackage{algorithmicx}
\usepackage[noend]{algpseudocode}
\usepackage{bm}
\begin{document}

\begin{frontmatter}

\title{On the linear convergence of distributed Nash equilibrium seeking for multi-cluster games under partial-decision information\thanksref{footnoteinfo}} 

\thanks[footnoteinfo]{This work was partially supported by Shanghai Pujiang Program under Grant 21PJ1413100, the National Natural Science Foundation of China under Grant 62003243, 62103305 and 62088101, by Shanghai Municipal Science and Technology Major Project under Grant 2021SHZDZX0100, Shanghai Municipal Commission of Science and Technology under Grant 19511132101, and Young Elite Scientist Sponsorship Program by cast of China Association for Science and Technology under Grant YESS20200136.}

\thanks[footnoteinfo]{Corresponding author: Xiuxian Li.}

\author[1,2]{Min Meng},
\ead{mengmin@tongji.edu.cn}
\author[1,2]{Xiuxian Li}
\ead{xli@tongji.edu.cn}
\address[1]{Department of Control Science and Engineering, College of Electronics and Information Engineering, Tongji University, Shanghai, China}
\address[2]{Shanghai Research Institute for Intelligent Autonomous Systems, Tongji University, Shanghai, China}

\begin{abstract}
This paper considers the distributed strategy design for Nash equilibrium (NE) seeking in multi-cluster games under a partial-decision information scenario. In the considered game, there are multiple clusters and each cluster consists of a group of agents. A cluster is viewed as a virtual noncooperative player that aims to minimize its local payoff function and the agents in a cluster are the actual players that cooperate within the cluster to optimize the payoff function of the cluster through communication via a connected graph. In our setting, agents have only partial-decision information, that is, they only know their own local cost functions, local feasible strategy sets and strategies of neighboring agents. To solve the NE seeking problem of this formulated game, a discrete-time distributed algorithm, called distributed projected gradient tracking algorithm (DPGT), is devised based on the inter- and intra-communication of clusters. In the designed algorithm, each agent is equipped with strategy variables including its own strategy and estimates of other clusters' strategies. With the help of a weighted Fronbenius norm and a weighted Euclidean norm, theoretical analysis is presented to rigorously show the linear convergence of the algorithm. Finally, a numerical example is given to illustrate the proposed algorithm.
\end{abstract}
\end{frontmatter}

\section{Introduction}
Game theory, which has been found to be a powerful tool to deal with optimization problems arising in multi-agent systems with the objective functions being coupled through decision variables of agents, has various applications including competitive markets \cite{li2015demand}, smart grids \cite{saad2012game}, transport systems \cite{hollander2006applicability}, to name just a few. A challenging issue in games is to design strategies to find a Nash equilibrium (NE) corresponding to the desirable and stable state, from which no agents want to deviate. Some references, such as \cite{facchinei2010generalized,yu2017distributed,shamma2005dynamic}, made an assumption that each agent can access all the competitors' decisions, which is impractical since a central node with bidirectional communication with all the players must exist in such case.

Therefore, in recent years, most scholars have focused on distributed NE seeking algorithms for noncooperative games composed of selfish decision-makers. For example, a payoff-based scheme was proposed in \cite{stankovic2012distributed,frihauf2011nash}, where each player is required to measure its cost function but not to communicate with others. In most circumstances, a player may not be aware of all other players' strategies, i.e., in a partial-decision information scenario. To handle such kind of partial-decision information scenarios, many results on the NE seeking problems were obtained both in continuous-time \cite{de2019distributed,gadjov2019passivity} and in discrete-time \cite{koshal2016distributed,salehisadaghiani2016distributed,salehisadaghiani2019distributed,tatarenko2019geometric}, in which gradient and consensus based algorithms were designed to estimate other players' strategies relying on local information. The algorithms in \cite{koshal2016distributed,salehisadaghiani2016distributed} equipped with vanishing stepsizes may have a slower convergence than those in \cite{salehisadaghiani2019distributed,tatarenko2019geometric,tatarenko2018geometric,bianchi2021fully} where fixed-step schemes were applied.

In contrast to noncooperative games, distributed optimization concerns a network of agents that collaborate to minimize the global cost function \cite{nedic2009distributed,zeng2017distributed,qu2018harnessing,li2020distri,li2019distributed_ijrnc,li2019distributed}. This problem is also an active research topic and has wide applications in resource allocation, machine learning, sensor networks, and energy systems \cite{bullo2009distributed}. Competition and cooperation among agents always coexist in many practical situations, such as health-care networks \cite{peng2009coexistence}, transportation networks \cite{shehory1998methods} and smart grids. For example, in smart grids, competition exists in energy management and energy market situations, while cooperation is reflected in economic dispatch. These practical situations may not be well modeled by only noncooperative games or distributed optimization problems. Inspired by the coexistence of competition and cooperation among agents, a multi-cluster (or multi-coalition) game was formulated in \cite{ye2018nash,ye2019unified,ye2020extremum,zeng2019generalized}. This game is conducted by multiple clusters (or coalitions), each of which is regarded as a virtual selfish player and aims to minimize its local payoff function. The agents in the same cluster are the actual players that cooperate within the cluster to optimize the payoff function of the cluster through communication via a connected graph.

However, all the above existing distributed NE seeking algorithms for multi-cluster games are in continuous-time. As discrete-time algorithms are easily implemented in practical applications, in this paper, we aim to design a distributed discrete-time algorithm for seeking an NE of multi-cluster games under partial-decision information. In the studied multi-cluster game, the payoff function of a cluster is defined as the average sum of local payoff functions of its agents and every cluster designates a representative agent to interact with other representative agents from other clusters through an arbitrary connected network.  With the aid of the available local information, each agent makes estimations of other clusters' strategies and the gradient of its cluster's payoff function at every iteration. Based on the inter- and intra-communication, a distributed projected gradient tracking algorithm (DPGT) is devised to find the NE of the studied multi-cluster game with strategy set constraints. Under some mild conditions, by introducing a weighted Frobenius norm and a weighted Euclidean norm, the algorithm is rigorously proved to converge to the NE at a linear rate. Finally, we present a numerical example of Cournot Competition games to illustrate the developed algorithm.

The main contributions of this paper can be summarized as follows:
\begin{itemize}
\item[1)] This paper is the first to design a discrete-time algorithm, in contrast to continuous-time algorithms in the existing references \cite{ye2018nash,ye2019unified,ye2020extremum,zeng2019generalized} for seeking NE of this kind of multi-cluster games with both competitive and cooperative behaviors. The challenge in devising the required algorithm is on the elaborate design of the weights of the information received from inter- and intra-communication.
    Additionally, the constraints on the strategy sets are also taken into account. In this setting, a novel gradient tracking and projection-based algorithm is proposed, which includes sonme existing algorithms for noncooperative games \cite{tatarenko2018geometric} and distributed optimization \cite{qu2018harnessing} as special cases.
\item[2)] The designed algorithm is proved to be not only convergent but also at a linear convergence rate, which is a significant challenge especially for the constrained multi-cluster games. It should be noted that the weighted matrix $\mathcal{A}$ corresponding to the connected communication graph among all agents is row-stochastic but unbalanced. Another novelty is that the left eigenvector $\pi$ with respect to eigenvalue 1 of $\mathcal{A}$ is explicitly computed, which plays an important role in dealing with the unbalanced communication and proving the linear convergence of the designed algorithm.
\end{itemize}

The rest of this paper is organized as follows. In Section \ref{section2}, the problem formulation is introduced.  The designed algorithm and main convergence result are presented in Section \ref{section3}.
Section \ref{section5} uses a numerical example to show the effectiveness of the proposed algorithm.  Section \ref{section6} makes a brief conclusion.

{\bf Notations}. Let $\mathbb{R}$, $\mathbb{R}^n$ and $\mathbb{R}^{m\times n}$ be the sets of real numbers, $n$-dimensional real column vectors and $m\times n$ real matrices, respectively. For an integer $n>0$, denote $[n]:=\{1,2,\ldots,n\}$. $I_n$ is the identity matrix of dimension $n$. ${\bf 1}_n$ (resp. ${\bf 0}_n$) represents an $n$-dimensional vector with all of its elements being 1 (resp. 0). For a vector or matrix $A$, $A^\top$ denotes the transpose of $A$ and $Col_i(A)$ (respectively $Row_i(A)$) is the $i$th column (respectively row) of $A$. $\rho(A)$ represents the spectral radius of $A$ and $det(A)$ is the determinant of $A$. For real symmetric matrices $P$ and $Q$, $P\succ(\succeq, \prec,\preceq)~Q$ means that $P-Q$ is positive (positive semi-, negative, negative semi-) definite, while for two vectors/matrices $w,v$ of the same dimension, $w\leq v$ means that each entry of $w$ is no greater than the corresponding one of $v$. $A\otimes B$ denotes the Kronecker product of matrices $A$ and $B$. ${\rm diag}\{a_1,a_2,\ldots,a_n\}$ or ${\rm diag}\{a_i,i\in[n]\}$ represents a diagonal matrix with $a_i$, $i\in[n]$, on its diagonal. For a vector $v$, we use ${\rm diag}(v)$ to represent the diagonal matrix with the vector $v$ on its diagonal. Denote by $col(z_1,\ldots,z_n)$ the column vector or matrix by pilling up $z_i$, $i\in[n]$.
Let $\langle w,v\rangle:=\sqrt{{\rm trace}[w^{\top}v]}$ for any two matrices or vectors $w,v$ of the same dimensions.
Denote by $\|\cdot\|_{F}$ and $\|\cdot\|$ the Frobenius norm and the Euclidean norm induced by the Frobenius inner product and the standard dot product, respectively, i.e., $\|v\|_{F}:=\sqrt{{\rm trace}[v^{\top}v]}$ and $\|w\|:=\sqrt{w^{\top}w}$ for $v\in\mathbb{R}^{n\times q}$ (or $\mathbb{R}^{n_i\times q_i}$) and $w\in\mathbb{R}^n$. The projection operator onto a set $\Omega$ is denoted by $P_{\Omega}[\cdot]$.

 An undirected graph, denoted as ${\mathcal G}=({\mathcal V},{\mathcal E},A)$, where ${\mathcal V}=\{1,2,\ldots,N\}$, ${\mathcal E}\subseteq{\mathcal V}\times{\mathcal V}$ and $A=(a_{ij})\in\mathbb{R}^{N\times N}$ represent the vertex set, the edge set and the weighted adjacency matrix of ${\mathcal G}$, respectively. The weights are defined as $a_{ij}>0$ if $(j,i)\in{\mathcal E}$ and $a_{ij}=0$ otherwise. $a_{ii}>0$ for all $i\in[N]$ in this paper.
 A path from node $i_1$ to node $i_l$ is composed of a sequence of edges $(i_h,i_{h+1})$, $h=1,2,\ldots,l-1$. An undirected graph ${\mathcal G}$ is said to be connected if for any vertices $i,j$, there is a path from node $i$ to node $j$.

\section{Problem formulation}\label{section2}
This paper is concerned with the multi-cluster noncooperative game, which is conducted by $m$ clusters. Each cluster $i\in[m]$ is a virtual self-interested player and contains $n_i$ ($\geq1$) agents communicating via an undirected graph $\mathcal{G}_i=([n_i],\mathcal{E}_i,A_i)$. In the meantime, each cluster designates a representative agent to interact with the representative agents from other clusters through an undirected communication topology $\mathcal{G}_0=([m],\mathcal{E}_0,A_0)$. Without loss of generality, it is supposed that the representative agent in cluster $i$ is the first agent, i.e., agent 1 in cluster $i$, where $i\in[m]$. The number of the agents in this game is $n:=\sum_{i=1}^mn_i$. The concepts of strategy variables and payoff functions of the multi-cluster game are given as follows.
\begin{itemize}
\item The {\em strategy variable} of agent $j$ in cluster $i$ is denoted as $x_{ij}\in\Omega_i\subset\mathbb{R}^{q_i}$, where $\Omega_i$ is non-empty, closed and convex. Let $x_i=col(x_{i1},\ldots,x_{in_i})$ be the strategy variable of cluster $i$ and $x_{-i}$ be the joint action of all the other clusters except that of $i$, i.e., $x_{-i}=col(x_1,\ldots,x_{i-1},x_{i+1},\ldots,x_m)$. The strategy variable of this game is defined as $x=col(x_1,\ldots,x_m)\in\mathbb{R}^N$, where $N:=\sum_{i=1}^mn_iq_i$.
\item The {\em payoff function} of cluster $i$, $f_i:\mathbb{R}^N\to\mathbb{R}$, is defined as
\begin{align}\label{equ1}
f_i(x_i,x_{-i})=\frac{1}{n_i}\sum\limits_{j=1}^{n_i}f_{ij}(x_{ij},\Gamma_{i}(x_{-i})),
\end{align}
where $f_{ij}(x_{ij},\Gamma_{i}(x_{-i})):\mathbb{R}^q\to\mathbb{R}$ is only available to agent $j$ in cluster $i$ and $\Gamma_{i}(x_{-i})\in\mathbb{R}^{q-q_i}$ is the stacked strategies of the representative agents of all the clusters except that of cluster $i$, i.e., $\Gamma_{i}(x_{-i})=col(x_{11},\ldots,x_{i-1,1},x_{i+1,1},\ldots,x_{m1})$. Here, $q:=\sum_{i=1}^mq_i$. Cluster $i\in[m]$ aims to seek a strategy $x_i^*=col(x_{i1}^*,\ldots,x_{in_i}^*)$ with $x_{ij}^*=x_{il}^*\in\Omega_i$ for $j,l\in[n_i]$ that minimizes its own payoff function $f_i(x_i,x_{-i})$ under $x_{-i}$.
\end{itemize}

Note that the strategies of agents in the same cluster are ensured to reach an agreement. A strategy profile $(x_i^*,x_{-i}^*)$ is called an NE of the formulated cluster game if $x_{ij}^*=x_{il}^*=\tilde{x}_i^*\in\Omega_i$ for all $j,l\in[n_i]$, and for all $i\in[m]$,
\begin{align}
f_i(x_i^*,x_{-i}^*)\leq f_i(x_i,x_{-i}^*),  ~x_i={\bf1}_{n_i}\otimes y_i,~ \forall y_i\in\Omega_i.
\end{align}
The objective of this paper is to design a distributed discrete-time algorithm to find an NE of the studied multi-cluster game based on local information though communication, i.e., under
a partial-decision information setting.

\begin{remark}
The formulated cluster game can model the coexistence of competition and cooperation simultaneously and subsume noncooperative games and distributed optimization as special cases. Specifically, if $n_i=1$ for all $i\in[m]$, the multi-cluster game is a noncooperative game among $m$ players \cite{koshal2016distributed,salehisadaghiani2016distributed,salehisadaghiani2019distributed,tatarenko2019geometric,tatarenko2018geometric}. If $m=1$, the considered problem is reduced to the distributed optimization problem, which has been investigated such as in \cite{nedic2009distributed,zeng2017distributed,qu2018harnessing,li2020distri,li2019distributed_ijrnc,li2019distributed}.
\end{remark}

To proceed, some standard assumptions are listed below.
\begin{assumption}\label{assumption1}
Graphs $\mathcal{G}_0,\mathcal{G}_1,\ldots,\mathcal{G}_m$ are undirected and connected. All adjacency matrices $A_0,A_1,\ldots,A_m$ are row and column stochastic, i.e., $A_0{\bf1}_m={\bf1}_m$, ${\bf1}_m^{\top}A_0={\bf1}_m^{\top}$, $A_i{\bf1}_{n_i}={\bf1}_{n_i}$, ${\bf1}_{n_i}^{\top}A_i={\bf1}_{n_i}^{\top}$, $i\in[m]$.
\end{assumption}
Under Assumption \ref{assumption1}, one has \cite{li2020}
\begin{align}\label{e3}
\sigma_i:=\|A_i-{\bf1}_{n_i}{\bf1}_{n_i}^{\top}/{n_i}\|<1.
\end{align}
Denote $\sigma_{\max}:=\max_{i\in[m]}\sigma_i$.
\begin{assumption}\label{assumption2}
For every $j\in[n_i]$, $i\in[m]$, local payoff function $f_{ij}(x_{ij},\Gamma_{i}(x_{-i}))$ is continuously differentiable and the gradient $\nabla_i f_{ij}(x_{ij},\Gamma_{i}(x_{-i})):=\frac{\partial f_{ij}(x_{ij},\Gamma_{i}(x_{-i}))}{\partial x_{ij}}$ is Lipschitz continuous on $\mathbb{R}^{q}$, i.e., for some constant $L_{ij}>0$,
\begin{align}
&\|\nabla_i f_{ij}(x_{ij},\Gamma_{i}(x_{-i}))-\nabla_i f_{ij}(\tilde{x}_{ij},\Gamma_{i}(\tilde{x}_{-i}))\|\nonumber\\
&\leq L_{ij}\left\|\left[\begin{array}{c}x_{ij}-\tilde{x}_{ij}\\ \Gamma_{i}(x_{-i}))-\Gamma_{i}(\tilde{x}_{-i})\end{array}\right]\right\|.
\end{align}
Denote $L:=\max_{j\in[n_i],i\in[m]}L_{ij}$.
\end{assumption}
Assumptions \ref{assumption1} and \ref{assumption2} are standard and commonly used in distributed discrete-time algorithms, such as distributed optimization, consensus, and NE seeking in noncooperative games \cite{koshal2016distributed,tatarenko2019geometric,tatarenko2018geometric,qu2018harnessing}.

Define
$
g_i(y_i,y_{-i}):=\frac{1}{n_i}\sum_{j=1}^{n_i}\nabla_if_{ij}(y_{i},y_{-i}),
$
where $y_i\in\mathbb{R}^{q_i}$, $i\in[m]$.
The game mapping $F:\mathbb{R}^q\to\mathbb{R}^q$ is defined as
\begin{align}\label{equ5}
F(y):=col(g_1(y_1,y_{-1}),\ldots,g_m(y_m,y_{-m})),~y\in\mathbb{R}^q.
\end{align}
\begin{assumption}\label{assumption3}
The mapping $F(y)$ is strongly monotone on $\mathbb{R}^q$ with constant $\mu>0$, that is, for any $y=col(y_1,\ldots,y_m)$ and $z=col(z_1,\ldots,z_m)$ with $y_i,z_i\in\mathbb{R}^{q_i}, i\in[m]$, the following inequality holds:
\begin{align}
&\sum\limits_{i=1}^m\frac{1}{n_i}\sum\limits_{j=1}^{n_i}\left[\nabla_if_{ij}(y_i,y_{-i})-\nabla_if_{ij}(z_i,z_{-i})\right]^{\top}(y_i-z_i)\nonumber\\
&\geq\mu\|y-z\|^2.
\end{align}
\end{assumption}
This assumption is standard in the literature on NE seeking algorithms with fast convergence rates and can ensure the existence and uniqueness of NE of the studied multi-cluster game, which is consistent with the condition ensuring the existence and uniqueness of NE of conventional noncooperative games in \cite{tatarenko2019geometric,tatarenko2018geometric}.

Note that in the multi-cluster game, the agents in the same cluster are required to reach an agreement to minimize the payoff function of the cluster, thus $x^*=(x_i^*,x_{-i}^*)$ is an NE, where $x_i^*=\mathbf{1}_{n_i}\otimes\tilde{x}_i^*$ and $\tilde{x}_i^*\in\Omega_i$, if and only if $\tilde{x}_i^*$ for any $i\in[m]$ satisfies
\begin{align}\label{equ9}
\Big\langle\frac{1}{n_i}\sum_{j=1}^{n_i}\nabla_if_{ij}(\tilde{x}_i^*,\tilde{x}_{-i}^*),x_i-\tilde{x}_i^*\Big\rangle\geq0,~\forall x_i\in\Omega_i,
\end{align}
which is equivalent to for any $\alpha_i>0$, it holds
\begin{align}\label{equ9e}
\tilde{x}_i^*
=P_{\Omega_i}\Big[\tilde{x}_i^*-\frac{\alpha_i}{n_i}\sum_{j=1}^{n_i}\nabla_if_{ij}(\tilde{x}_i^*,\tilde{x}_{-i}^*)\Big], ~i\in[m].
\end{align}

\section{Main result}\label{section3}
In this section, a distributed discrete-time algorithm is proposed for NE seeking of the considered multi-cluster game. To deal with the partial-decision information scenario, where the agents can access the local information exchanged through local communication, it is assumed that each agent $j$ in cluster $i$ maintains a local variable
\begin{align}
x_{(ij)}=col(x_{(1)ij},\ldots,x_{(m)ij})\in\mathbb{R}^q,
\end{align}
which is the estimation of the joint strategy of the representative agents, i.e., $col(x_{11},\ldots,x_{m1})$. Here $q=\sum_{i=1}^mq_i$, $x_{(s)ij}\in\mathbb{R}^{q_s}$ is the estimate of $x_{s1}$ at the agent $j$ in cluster $i$ and $x_{(i)ij}=x_{ij}$. Also, the estimates of other representative agents except that of cluster $i$ is compactly denoted by
\begin{align}
x_{-(ij)}=col(x_{(1)ij},\ldots,x_{(i-1)ij},x_{(i+1)ij},\ldots,x_{(m)ij}).
\end{align}
A distributed projected gradient tracking algorithm (DPGT) is proposed as in Algorithm \ref{alg1} to learn the NE of the multi-cluster game.
\begin{algorithm}[!htbp]\caption{Distributed Projected Gradient Tracking Algorithm (DPGT)}\label{alg1}
Each agent $j$ in cluster $i$ maintains vector variables $x_{(ij)}^t=col(x^t_{(1)ij},\ldots,x_{(m)ij}^t)\in\mathbb{R}^q$ and $v_{ij}^t\in\mathbb{R}^{q_i}$ at iteration $t$.

 {\bf Initialization:} Initialize $x_{ij}^0\in\Omega_i$, $x_{-(ij)}^0\in\mathbb{R}^{q-q_i}$ arbitrarily and let $v_{ij}^0=\nabla_if_{ij}(x_{ij}^0,x_{-(ij)}^0)$.

{\bf Iteration:} For $t\geq 0$, every agent $j$ in cluster $i$ processes the following update:
\begin{subequations}\label{equ10}
\begin{align}
x_{ij}^{t+1}&=\left\{\begin{array}{ll}P_{\Omega_i}\Big[\frac{1}{2}\sum\limits_{l=1}^{n_i}a_i^{jl}x_{il}^t
+\frac{1}{2}\sum\limits_{h=1}^ma_0^{ih}x_{(i)hj}^t-\alpha_iv_{ij}^t\Big],&j=1,\\
P_{\Omega_i}\Big[\sum\limits_{l=1}^{n_i}a_i^{jl}x_{il}^t-\alpha_iv_{ij}^t\Big],&j\neq1,
\end{array}\right.\label{}\\
x_{(s)ij}^{t+1}&=\left\{\begin{array}{ll}
\frac{1}{2}\sum\limits_{l=1}^{n_i}a_i^{jl}x_{(s)il}^t+\frac{1}{2}\sum\limits_{h=1}^ma_0^{ih}x_{(s)hj}^t,&j=1,s\neq i,\\
\sum\limits_{l=1}^{n_i}a_i^{jl}x_{(s)il}^t,&j\neq1,s\neq i,
\end{array}\right.\label{}\\
v_{ij}^{t+1}&=\sum\limits_{l=1}^{n_i}a_i^{jl}v_{il}^t+\nabla_if_{ij}(x_{ij}^{t+1},x_{-(ij)}^{t+1})-\nabla_if_{ij}(x_{ij}^{t},x_{-(ij)}^{t}),\label{equ10c}
\end{align}
\end{subequations}
where $a_0^{ih}$ is the $(i,h)$ element of $A_0$, $a_i^{jl}$ is the $(j,l)$ element of $A_i$, $i\in[m]$, and $\alpha_i>0$, $i\in[m]$, are the stepsizes to be determined.
\end{algorithm}

In Algorithm \ref{alg1}, each agent $j$ in cluster $i$ can only access its own strategy set $\Omega_i$, its own strategy variable $x_{ij}^t$, the estimation $x_{-(ij)}^t$ of strategies of representative agents from other clusters, the value $\nabla_{ij}f_{ij}(x_{ij}^t,x_{-(ij)}^t)$, the gradient estimation variable $v_{ij}^t$, and the neighbors' strategies from inter/intra-communication.
Algorithm \ref{alg1} is designed relying on three parts: inter-cluster update mechanism, intra-cluster update mechanism and estimation of the gradients of local clusters' payoff functions. In Algorithm \ref{alg1}, the inter-cluster update mechanism is reflected at variables $x_{(i1)}^t$, $i\in[m]$, since only the representative agents of clusters can communicate with each other via communication topology $\mathcal{G}_0$. Meanwhile, the update of $x_{(i1)}^t$, $i\in[m]$ should also combine the intra-cluster communication via graph $\mathcal{G}_i$, $i\in[m]$. The non-representative agents only need to update their local variables by considering intra-cluster communication. That is, $x_{(ij)}^t$ for $j\neq1$, $i\in[m]$ follow the intra-cluster update mechanism. $v_{ij}^t$ is an auxiliary variable to estimate the gradient of the payoff function of cluster $i$ at the estimated strategy of other representative agents, i.e., $\frac{1}{n_i}\sum_{j=1}^{n_i}\nabla_if_{ij}(x_{ij}^t,x_{-(ij)}^t)$.
Every agent updates its local variables $x_{(ij)}^t$ and $v_{ij}^t$ only by local information. Thus this algorithm is distributed.

For convenient analysis, let us rewrite (\ref{equ10}) into a compact form in the following. Define the estimation matrix ${\bf x}^t\in\mathbb{R}^{n\times q}$ with $n:=\sum_{i=1}^mn_i$ as
\begin{align}
{\bf x}^t:=col({\bf x}_1^t,\ldots,{\bf x}_m^t),\label{e14}
\end{align}
where
$
{\bf x}_i^t:=[\begin{array}{cccc}x_{(i1)}^t&x_{(i2)}^t&\cdots&x_{(in_i)}^t\end{array}]^{\top}\in\mathbb{R}^{n_i\times q},~i\in[m].
$
Similarly, denote
\begin{align*}
{\bf v}_i^t&:=[\begin{array}{cccc}v_{i1}^t&v_{i2}^t&\ldots&v_{in_i}^t\end{array}]^{\top}\in\mathbb{R}^{n_i\times q_i},~i\in[m],\\
 V^t&:={\rm diag}\{{\bf v}_1^t,\ldots,{\bf v}_m^t\}\in\mathbb{R}^{n\times q},\\
 \Lambda&:={\rm diag}\{\alpha_1I_{n_1},\ldots,\alpha_mI_{n_m}\}.
\end{align*}
Considering all the $n$ agents in the multi-cluster game, the entire communication topology $\mathcal{G}$ among the $n$ agents is composed of graphs $\mathcal{G}_0,\mathcal{G}_1,\ldots,\mathcal{G}_m$ with the vertex set being $[n]$ composed of all the $n$ agents and the adjacency matrix $\mathcal{A}\in\mathbb{R}^{n\times n}$ being given as
\begin{align}\label{equ13}
\mathcal{A}=\left[\begin{array}{ccccc}\tilde{A}_1+B_{11}&B_{12}&\cdots&B_{1m}\\B_{21}&\tilde{A}_2+B_{22}&\cdots&B_{2m}\\ \vdots&\vdots&&\vdots\\B_{m1}&B_{m2}&\cdots&\tilde{A}_m+B_{mm}\end{array}\right],
\end{align}
where $\tilde{A}_i\in\mathbb{R}^{n_i\times n_i}$ is equal to the adjacency matrix $A_i$ of $\mathcal{G}_i$ except that the first row $Row_1(\tilde{A}_i)$ is $\frac{1}{2}Row_1(A_i)$ and $B_{ij}\in\mathbb{R}^{n_i\times n_j}$ has the $(1,1)$ element as $\frac{1}{2}a_0^{ij}$ and other elements as $0$, $i,j\in[m]$. Then $\mathcal{A}$ is a row stochastic matrix and has positive diagonal entries. For a matrix $A\in\mathbb{R}^{n\times q}$, denote by $\mathcal{P}_{\Omega}[A]$ a matrix with its $h$th row being
\begin{align}
(P_{\bm{\Omega}_h}[Col_h(A^{\top})])^{\top},~h\in[n],\label{e19}
\end{align}
where $\bm{\Omega}_h:=\{y\in\mathbb{R}^q\mid Q_iy\in\Omega_i\}$ for $h\in[n_{<i}+1,n_{<i}+n_i]$ and $Q_i:=[\mathbf{0}_{n_i\times n_{<i}},I_{n_i},\mathbf{0}_{n_i\times n_{>i}}]$ with $n_{<i}:=\sum_{j=1}^{i-1}n_j$, $n_{>i}:=\sum_{j=i+1}^{m}n_j$, $i\in[m]$, and $P_{\bm{\Omega}_h}[\cdot]$ represents the projection operator onto the set $\bm{\Omega}_h$. By the above notations, (\ref{equ10}) can be rewritten into a compact form as follows:
\begin{align}
{\bf x}^{t+1}&=\mathcal{P}_{\Omega}\left[\mathcal{A}{\bf x}^t-\Lambda V^t\right],\label{equ14}\\
{\bf v}_i^{t+1}&=A_i{\bf v}_i^t+G_i({\bf x}_i^{t+1})-G_i({\bf x}_i^t),\label{equ15}
\end{align}
where for $i\in[m]$,
\begin{align}\label{e18}
G_i({\bf x}_i^t)&=[\begin{array}{cccc}\nabla_if_{i1}(x_{i1}^t,x_{-(i1)}^t)& \cdots & \nabla_if_{in_i}(x_{in_i}^t,x_{-(in_i)}^t)\end{array}]^{\top}.
\end{align}

It can be seen that $\mathcal{A}$ is a row stochastic matrix with a simple eigenvalue $1$ under Assumption \ref{assumption1}. Then, $\mathcal{A}$ has a left eigenvector $\pi\in\mathbb{R}^n$ corresponding to eigenvalue 1 satisfying that ${\pi}^{\top}\mathcal{A}=\pi^{\top}$, ${\pi}^{\top}{\bf1}_n=1$ and every element of $\pi$ is positive. Indeed, by the structure of $\mathcal{A}$ in (\ref{equ13}), $\pi$ can be obtained in the following lemma. 
\begin{lemma}\label{lem1}
Under Assumption \ref{assumption1}, the left eigenvector $\pi\in\mathbb{R}^n$ of matrix $\mathcal{A}$ in (\ref{equ13}) corresponding to eigenvalue 1 such that ${\pi}^{\top}\mathcal{A}=\pi^{\top}$ and ${\pi}^{\top}{\bf1}_n=1$ is given as $\pi=col(\pi^1,\ldots,\pi^m)$ with $\pi^i$ being $\pi^i=col(\frac{2}{n+m},\frac{1}{n+m},\ldots,\frac{1}{n+m})\in\mathbb{R}^{n_i}$, $i\in[m]$. Moreover, ${\bf 1}_{n_i}^{\top}\pi^i=\frac{n_i+1}{n+m}$, $i\in[m]$.
\end{lemma}
{\bf Proof.} The proof is postponed to Appendix \ref{appendix_A}. \hfill$\blacksquare$

Based on $\pi=col(\pi_1,\ldots,\pi_n)$, a weighted Euclidean norm and a weighted Frobenius norm are, respectively, defined as follows: for $x,y\in\mathbb{R}^n$ and ${\bf x},{\bf y}\in\mathbb{R}^{n\times q}$,
\begin{align}
\|x\|_{\pi}&:=\sqrt{\langle{\rm diag}({\pi}){x},{x}\rangle}=\|{\rm diag}(\sqrt{\pi}){x}\|,\label{equ18}\\
\|{\bf x}\|_{F}^{\pi}&:=\sqrt{\langle{\rm diag}({\pi}){\bf x},{\bf x}\rangle}=\|{\rm diag}(\sqrt{\pi}){\bf x}\|_{F},\label{equ19}
\end{align}
where $\sqrt{\pi}:=col(\sqrt{\pi}_1,\ldots,\sqrt{\pi}_n)$.
We also denote by $\|B\|_{\pi}$ the matrix norm of matrix $B\in\mathbb{R}^{n\times n}$ induced by the weighted Euclidean norm. Then
\begin{align}
\|B\|_{\pi}=\|{\rm diag}(\sqrt{\pi})B{\rm diag}(\sqrt{\pi})^{-1}\|.
\end{align}
The following results can be obtained.
\begin{lemma}\label{lemma1}
\begin{itemize}
\item[1)] For two positive semi-definite matrices $P,Q\in\mathbb{R}^{n\times n}$, if $P\succeq Q$, then ${\rm trace}[P]\geq{\rm trace}[Q]$. Moreover, for any $A\in\mathbb{R}^{n\times n}$, $B\in\mathbb{R}^{n\times q}$, we have
    $\|AB\|_{F}\leq\|A\|\|B\|_{F}$.
\item[2)] The Euclidean norm $\|\cdot\|$ and the weighted Euclidean norm $\|\cdot\|_{\pi}$ are equivalent. Specifically, for any $x\in\mathbb{R}^n$, there holds:
\begin{align}\label{equ21}
\frac{1}{\sqrt{n+m}}\|x\|\leq\|x\|_{\pi}\leq\frac{\sqrt{2}}{\sqrt{n+m}}\|x\|.
\end{align}
\item[3)] The Frobenius norm $\|\cdot\|_{F}$ and the weighted Frobenius norm $\|\cdot\|_{F}^{\pi}$ are equivalent. Specifically, for any ${\bf x}\in\mathbb{R}^{n\times q}$, the following inequality holds:
 \begin{align}\label{equ22}
 \frac{1}{\sqrt{n+m}}\|{\bf x}\|_{F}\leq\|{\bf x}\|_{F}^{\pi}\leq\frac{\sqrt{2}}{\sqrt{n+m}}\|{\bf x}\|_{F}.
 \end{align}
 \item[4)] Under Assumption \ref{assumption1}, one has
 \begin{align}
\|\mathcal{A}x-\mathcal{A}_{\infty}x\|_{{\pi}}&\leq\sigma\|x-\mathcal{A}_{\infty}x\|_{{\pi}},~\forall x\in\mathbb{R}^n,\label{equ23}\\
\|\mathcal{A}{\bf x}-\mathcal{A}_{\infty}{\bf x}\|_{F}^{{\pi}}&\leq\sigma\|{\bf x}-\mathcal{A}_{\infty}{\bf x}\|_{F}^{{\pi}},~\forall{\bf x}\in\mathbb{R}^{n\times q},\label{equ24}
 \end{align}
where $\mathcal{A}_{\infty}:={\bf1}_n\pi^{\top}$ and $\sigma:=\|\mathcal{A}-\mathcal{A}_{\infty}\|_{\pi}<1$.
\end{itemize}
\end{lemma}
{\bf Proof.} The proof is postponed to Appendix \ref{appendix_B}. \hfill$\blacksquare$

Note that it is required that all the agents in the same cluster take the same strategy when reaching the NE $x^*$, namely, the strategies $x_{ij}^*$ and $x_{il}^*$ of agent $j$ and agent $l$ in cluster $i$, respectively, should be equal. For ease of notations, we can denote the NE as $\tilde{x}^*=col(\tilde{x}_1^*,\ldots,\tilde{x}_m^*)\in\mathbb{R}^q$, where $q=\sum_{i=1}^mq_i$, and $\tilde{x}_i^*=x_{ij}^*$ for $j\in[n_i]$ and $i\in[m]$. Denote $\tilde{\bf x}^*:={\bf1}_n(\tilde{x}^*)^{\top}$.
 Then it is ready to present the main result on the convergence of Algorithm \ref{alg1}.
\begin{theorem}\label{thm1}
If Assumptions \ref{assumption1}--\ref{assumption3} hold and the stepsizes are chosen as $\alpha_i=\alpha/(n_i+1)$ for $\alpha>0$, $i\in[m]$, then $x_{ij}^t$ and $v_{ij}^t$ generated by Algorithm \ref{alg1} satisfy
\begin{align}
\left[
\begin{array}{c}
\|{\bf x}^{t+1}-\tilde{\bf x}^*\|_F^{\pi}\\
\sum\limits_{i=1}^m\|{\bf v}_i^{t+1}-\overline{\bf v}_i^{t+1}\|_F
\end{array}
\right]
\leq H_{\alpha}
\left[
\begin{array}{c}
\|{\bf x}^{t}-\tilde{\bf x}^*\|_F^{\pi}\\
\sum\limits_{i=1}^m\|{\bf v}_i^{t}-\overline{\bf v}_i^{t}\|_F
\end{array}
\right],
\end{align}
where $\overline{\bf v}^t_i:=\frac{{\bf1}_{n_i}{\bf 1}_{n_i}^{\top}}{n_i}{\bf v}_i^t\in\mathbb{R}^{n_i\times q_i}$, and
\begin{align}
&H_{\alpha}:=\left(
\begin{array}{cc}\sqrt{\rho(M_{\alpha})}&\frac{\sqrt{2}\alpha}{2\sqrt{n+m}}\\
\sqrt{m(n+m)}(1+\sqrt{\rho(M_{\alpha})})L&~\sigma_{\max}+\frac{\sqrt{2m}L\alpha}{2}
\end{array}
\right),\\
&M_{\alpha}:=\left(
\begin{array}{cc}1-\frac{2\mu\alpha}{n}+\frac{L^2\alpha^2}{2}&\frac{\sqrt{2}(1+\sigma)L\alpha}{2}\\
  \frac{\sqrt{2}(1+\sigma)L\alpha}{2}&\sigma^2+\sqrt{2}\sigma L\alpha+\frac{L^2\alpha^2}{2}
\end{array}
\right).
\end{align}
Moreover, if $0<\alpha<1$ is chosen such that $\rho(M_{\alpha})<1$ and $\rho(H_{\alpha})<1$, then $x_{ij}^t$ converges to the NE $x_{ij}^*$ linearly with the convergence rate being $\rho(H_{\alpha})$.
\end{theorem}
{\bf Proof.} The proof is postponed to Appendix \ref{appendix_C}. \hfill$\blacksquare$

It can be seen from Theorem \ref{thm1} that the sequence $\{x^t_{ij}\}$ generated by Algorithm \ref{alg1} converges to the NE $x_{ij}^*$ with a linear convergence rate $\rho(H_{\alpha})$ if the stepsizes are appropriately chosen. Also, the convergence rate depends on the Lipschitz constant $L$, the strongly monotonicity constant $\mu$, the number of agents, the dimension of strategies, and the properties of the communication topology. The following result present a sufficient condition on $\alpha$ ensuring $\rho(M_{\alpha})<1$ and $\rho(H_{\alpha})<1$.
\begin{proposition}\label{pro1}
For the matrices $M_{\alpha}$ and $H_{\alpha}$ defined in Theorem \ref{thm1}, $\rho(M_{\alpha})<1$ and $\rho(H_{\alpha})<1$ hold if $\alpha$ satisfies
\begin{align}
0<\alpha&<\min\left\{\frac{n}{4\mu},\frac{1}{\sqrt{2}(1+\sigma)L},\frac{\sqrt{2}\sigma^2}{(1-\sigma)L},\right.\nonumber\\
&\left.~~~\frac{1-\sigma^2}{3\sqrt{2}\sigma L},\frac{2\mu(1-\sigma)}{3nL^2},\frac{\sqrt{2(1-\sigma^2)}}{\sqrt{3}L},\right.\nonumber\\
&~~~\left.\frac{1-\sigma_{\max}}{\sqrt{2m}L},\frac{(1-\sigma_{\max})(1-\sqrt{\rho(M_{\alpha})})}{\sqrt{2m}(1+\sqrt{\rho(M_{\alpha})})L}\right\}.\label{e32}
\end{align}
\end{proposition}
{\bf Proof.} The proof is postponed to Appendix \ref{appendix_D}. \hfill$\blacksquare$

A sufficiently small $\alpha$ satisfying (\ref{e32}) can be always found. The bounds in (\ref{e32}) on $\alpha$ are not tight, and in practice better bounds on $\alpha$ can be obtained by simply verifying whether $\rho(M_{\alpha})<1$ and $\rho(H_{\alpha})<1$ are satisfied.

\begin{remark}
If there is only one cluster in the studied multi-cluster game, i.e., $m=1$, then Algorithm \ref{alg1} and Theorem \ref{thm1} will reduce to the result for distributed optimization, which is consistent with that in \cite{qu2018harnessing}. If there is only one agent in each cluster, i.e., $n_i=1$, $i\in[m]$, then the studied game becomes a conventional noncooperative game in \cite{koshal2016distributed,salehisadaghiani2016distributed,salehisadaghiani2019distributed,tatarenko2019geometric} and Algorithm \ref{alg1} is the gradient-based algorithm in \cite{tatarenko2018geometric}.
\end{remark}
\begin{remark}
In comparison, the existing NE seeking algorithms for multi-cluster games in \cite{ye2018nash,ye2019unified,ye2020extremum,zeng2019generalized} were designed in continuous-time and under full-decision information, while Algorithm \ref{alg1} here is a discrete-time algorithm under a partial-decision information setting, and meanwhile, explicit upper bounds on $\alpha$ are provided to ensure a linear convergence rate. Although discrete-time algorithms may be obtained by discretizing continuous-time algorithms using such as explicit Euler method, it is unclear how to select the sampling stepsize to guarantee the convergence, especially when ensuring the linear convergence.
\end{remark}

\section{Example}\label{section5}
In this section, we present a numerical example to illustrate our algorithm. To this end, a Cournot competition game is considered as follows.

There are 5 father companies, which are regarded as clusters in the multi-cluster game, and each father company has 20 subsidiary companies, which are viewed as agents in clusters. The father companies compete with each other by adjusting the production quantity of goods. The subsidiary companies affiliated with the same father company produce components for this father company, cooperate to reach an agreement and meanwhile ensure that the profit of the father company is optimal. Assume that the subsidiary companies in the same father company $i$ can communicate through a connected graph $\mathcal{G}_i$ and each father company appoints a subsidiary company to contact with other representative subsidiary companies from other father companies through another connected graph $\mathcal{G}_0$. By regarding the quantities of goods as strategy variables, denote by $x_{ij}$ the quantity of goods of subsidiary company $j$ in father company $i$, which is in $[0,30]$. The cost for producing components of goods and the price of components for per unit product by subsidiary company $j$ in father company $i$ are assumed to be $c_{ij}=5x_{ij}^2+i5x_{ij}+i$ and $p_{ij}=60i-\sum_{h=1}^5a_0^{ih}x_{h1}$, respectively, where $A_0=(a_0^{ih})\in\mathbb{R}^{5\times5}$ is the adjacency matrix of graph $\mathcal{G}_0$, $j\in\{1,2,\ldots,20\}$, $i,h\in\{1,2,3,4,5\}$. Then the payoff function of father company $i$ is $f_i=\sum_{j=1}^{20}f_{ij}$, where $f_{ij}=c_{ij}-x_{ij}p_{ij}$.

\begin{figure}[!ht]
 \centering
  \includegraphics[width=3in]{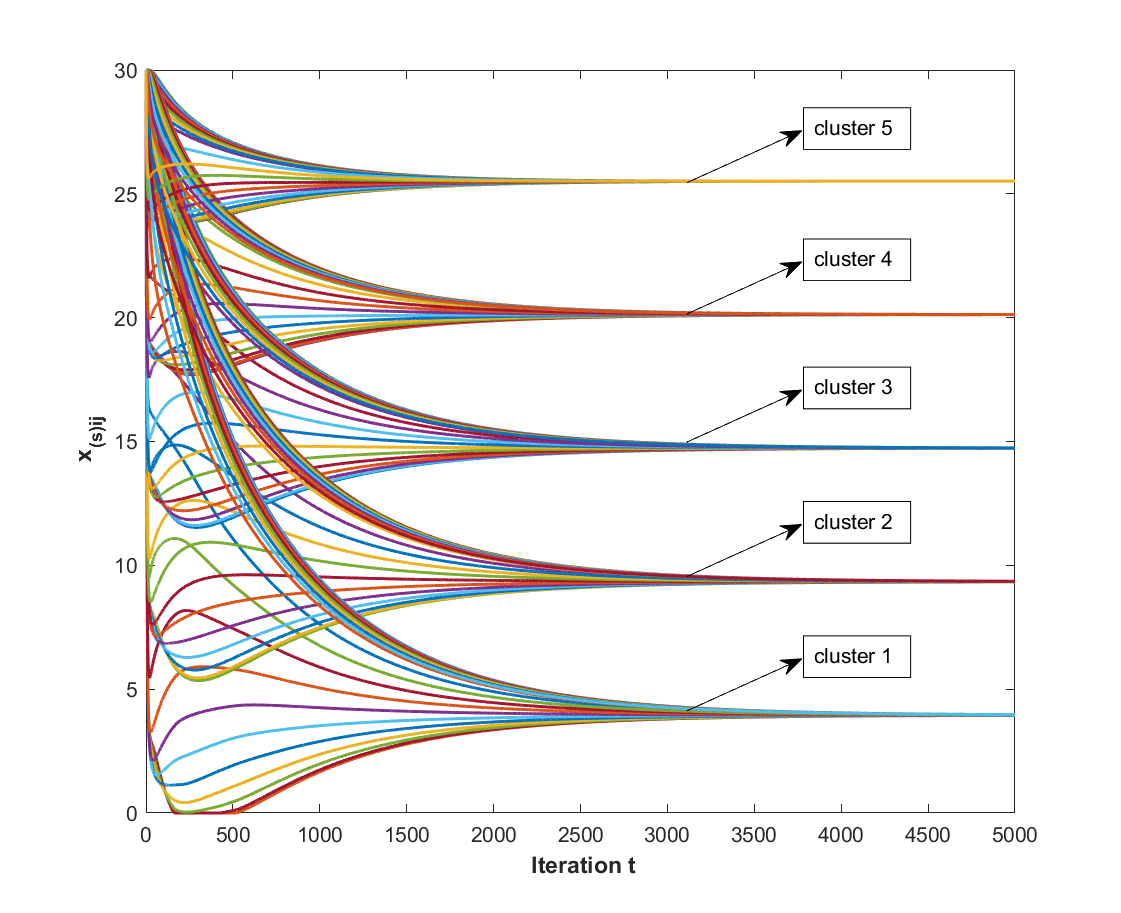}
  \caption{Trajectories of $x_{(s)ij}^t$, $s,i\in\{1,2,3,4,5\}$, $j\in\{1,2,\ldots,20\}$.}\label{fig1}
\end{figure}

This can be modeled as a multi-cluster game studied in this paper. It can be verified that Assumptions \ref{assumption1}--\ref{assumption3} are satisfied. By (\ref{equ9}) and a centralized method, the unique NE can be calculated as
$x_{1j}^*=3.9478$, $x_{2j}^*=9.3400$, $x_{3j}^*=14.7321$, $x_{4j}^*=20.1243$, $x_{5j}^*=25.5165$, where $j\in\{1,2,\ldots,20\}$. By our proposed algorithm with $\alpha_i=\alpha/(n_i+1)$ and $\alpha=0.2$, the sequence $\{x_{ij}^t\}$, as well as the estimate sequence $\{x_{(s)ij}^t\}$ of the strategy of cluster $s$ by agent $j$ in cluster $i$, converges to the unique NE, as shown in Fig. \ref{fig1}.

\begin{figure}[!ht]
 \centering
  \includegraphics[width=3in]{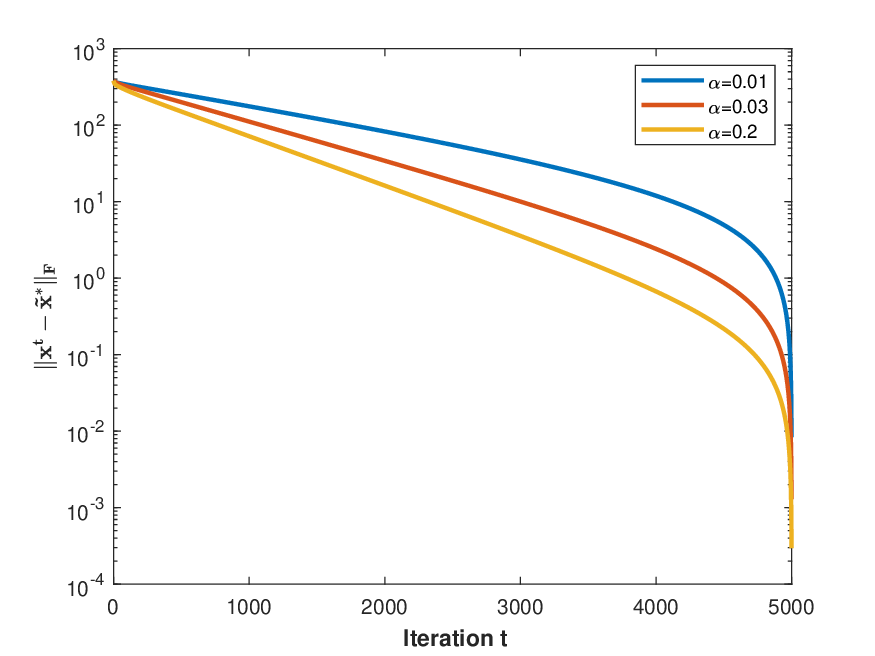}
  \caption{Influence of different $\alpha$ on the convergence rate.}\label{fig2}
\end{figure}

\begin{figure}[!ht]
 \centering
  \includegraphics[width=3in]{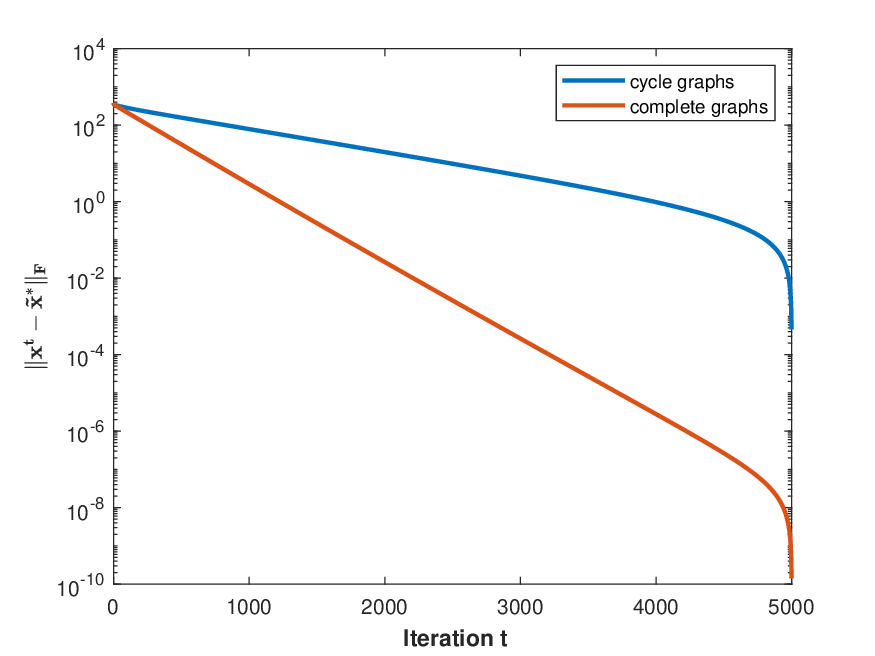}
  \caption{Influence of communication graphs on the convergence rate. The blue and red lines are for the cases, where inter- and intra- communication graphs are cycle and complete graphs, respectively.}\label{fig3}
\end{figure}

Furthermore, the evolutions of $\|{\bf x}^t-\tilde{\bf x}^*\|_F$ are shown in Fig. \ref{fig2} and Fig. \ref{fig3} under different stepsizes and different communication graphs, respectively, from which one can see that the convergence is linear, and larger stepsize and denser communication graphs lead to faster convergence rates. 
\section{Conclusion}\label{section6}
In this paper, we have studied the distributed NE seeking problem for a class of multi-cluster noncooperative games under a partial-decision information scenario. To design a distributed algorithm, every agent needs to make estimations of the strategies of other clusters at each iteration since each agent only has access to its local strategy set, its local payoff function coupled with other clusters' strategies and the neighbors' information. Then based on the inter- and intra-communication of clusters, a distributed projected gradient tracking algorithm in discrete-time was devised to find the unique NE of the multi-cluster game. Rigorous convergence analysis with a linear convergence rate was provided by introducing a weighted Frobenius norm and a weighted Euclidean norm. To further study generalized NE seeking for the formulated multi-cluster noncooperative games with constrained action sets and inequality constraints under a partial-decision information scenario is an interesting future research direction.
\section*{Appendix}
\appendix
\section{Proof of Lemma \ref{lem1}:}\label{appendix_A}
By splitting the adjacency matrices $A_i$ of $\mathcal{G}_i$ as
\begin{align}
A_i=\left[\begin{array}{cccc}A^i_{11}&A^i_{12}\\A^i_{21}&A^i_{22}\end{array}\right],
\end{align}
where $A^i_{11}\in\mathbb{R}$, $A^i_{12}\in\mathbb{R}^{1\times(n_i-1)}$, $A^i_{21}\in\mathbb{R}^{n_i-1}$, $A^i_{22}\in\mathbb{R}^{(n_i-1)\times(n_i-1)}$ and $i\in[m]$, one can rewrite $\tilde{A}_i$ as
\begin{align}
\tilde{A}_i=\left[\begin{array}{cccc}\frac{1}{2}A^i_{11}&\frac{1}{2}A^i_{12}\\A^i_{21}&A^i_{22}\end{array}\right], ~i\in[m].
\end{align}
Accordingly, $B_{hi}$ can be partitioned as
\begin{align}
B_{hi}=\left[\begin{array}{cccc}\frac{1}{2}a_0^{hi}&{\bf0}_{1\times(n_i-1)}\\{\bf0}_{n_h-1}&{\bf0}_{(n_h-1)\times(n_i-1)}\end{array}\right], ~h,i\in[m].
\end{align}
Let $\pi^i=(\pi^i_1,(\pi^i_2)^{\top})^{\top}$, where $\pi^i_1\in\mathbb{R}$, $i\in[m]$. Then from $\pi^{\top}\mathcal{A}=\pi^{\top}$, it can be obtained that for $i\in[m]$,
\begin{align*}
(\pi^i)^{\top}\tilde{A}_i+\sum\limits_{h=1}^m(\pi^h)^{\top}B_{hi}=(\pi^i)^{\top},
\end{align*}
which is equivalent to
\begin{align}
&[\pi^i_1~~(\pi^i_2)^{\top}]\left[\begin{array}{cccc}\frac{1}{2}A^i_{11}&\frac{1}{2}A^i_{12}\\A^i_{21}&A^i_{22}\end{array}\right]
+\sum\limits_{h=1}^m[\pi^h_1~~(\pi^h_2)^{\top}]\left[\begin{array}{cccc}\frac{1}{2}a_0^{hi}&{\bf0}\\{\bf0}&{\bf0}\end{array}\right]
\nonumber\\
&=[\pi^i_1~~(\pi^i_2)^{\top}].\label{e22}
\end{align}
Note that $A_i{\bf1}_{n_i}={\bf1}_{n_i}$, then right multiplying ${\bf1}_{n_i}$ on both sides of (\ref{e22}) yields that
\begin{align*}
&[\pi^i_1~~(\pi^i_2)^{\top}]\left[\begin{array}{cccc}\frac{1}{2}\\{\bf1}_{n_i-1}\end{array}\right]
+\sum\limits_{h=1}^m[\pi^h_1~~(\pi^h_2)^{\top}]\left[\begin{array}{cccc}\frac{1}{2}a_0^{hi}\\{\bf0}_{n_i-1}\end{array}\right]
\nonumber\\
&=\pi^i_1+(\pi^i_2)^{\top}{\bf1}_{n_i-1},
\end{align*}
that is,
\begin{align*}
\frac{1}{2}\pi^i_1+(\pi^i_2)^{\top}{\bf1}_{n_i-1}+\frac{1}{2}\sum\limits_{h=1}^m\pi^h_1a_0^{hi}=\pi^i_1+(\pi^i_2)^{\top}{\bf1}_{n_i-1},~i\in[m].
\end{align*}
Thus, $\sum\limits_{h=1}^m\pi^h_1a_0^{hi}=\pi^i_1,~i\in[m],$ from which one can see
 \begin{align}
(\pi^1_1,\ldots,\pi^m_1)A_0=(\pi^1_1,\ldots,\pi^m_1).
\end{align}
Therefore, the vector $col(\pi^1_1,\ldots,\pi^m_1)$ constructed by the first elements of $\pi^i$, $i\in[m]$, is a left eigenvector of $A_0$ corresponding to eigenvalue 1. Under Assumption \ref{assumption1}, since $\mathcal{G}_0$ is a connected graph and $A_0$ is a row and column stochastic matrix, $col(\pi^1_1,\ldots,\pi^m_1)$ should satisfy that $\pi^1_1=\cdots=\pi^m_1=\varpi$ for some $\varpi\in\mathbb{R}$. Consequently, (\ref{e22}) can be rewritten as
\begin{align}
&\left[\frac{1}{2}\varpi~~(\pi^i_2)^{\top}\right]\left[\begin{array}{cccc}A^i_{11}&A^i_{12}\\A^i_{21}&A^i_{22}\end{array}\right]
+\sum\limits_{h=1}^m\left[\frac{1}{2}\varpi a_0^{hi}~~{\bf0}_{1\times(n_i-1)}\right]
\nonumber\\
&=[\varpi~~(\pi^i_2)^{\top}].\label{e23}
\end{align}
Also from Assumption \ref{assumption1}, ${\bf1}_m^{\top}A_0={\bf1}_m^{\top}$, then for $i\in[m]$, $\sum_{h=1}^ma_0^{hi}=1$ and (\ref{e23}) becomes
\begin{align*}
\left[\frac{1}{2}\varpi~~(\pi^i_2)^{\top}\right]A_i+\left[\frac{1}{2}\varpi~~{\bf0}_{n_i-1}\right]
=\left[\varpi~~(\pi^i_2)^{\top}\right],
\end{align*}
that is,
\begin{align}
\left[\frac{1}{2}\varpi~~(\pi^i_2)^{\top}\right]A_i
=\left[\frac{1}{2}\varpi~~(\pi^i_2)^{\top}\right].
\end{align}
Hence, for $i\in[m]$, $\left(\frac{1}{2}\varpi,(\pi^i_2)^{\top}\right)$ is a left eigenvector of $A_i$ corresponding to eigenvalue 1, thus, under Assumption \ref{assumption1}, one has
$
\pi^i_2=\frac{1}{2}\varpi{\bf1}_{n_i-1}, ~i\in[m].
$
Note that $1=\pi^{\top}{\bf1}_n$, then it is obtained that
$
1=\sum_{i=1}^m(1+\frac{n_i-1}{2})\varpi=\frac{n+m}{2}\varpi,
$
indicating $\varpi=\frac{2}{n+m}$ and $\pi^i=col(\frac{2}{n+m},\frac{1}{n+m}\ldots,\frac{1}{n+m})$. The last claim can be easily verified and the proof is thus completed. \hfill$\blacksquare$

\section{Proof of Lemma \ref{lemma1}:}\label{appendix_B}
{\bf Proof.} 1)--3) can be easily proved by the definitions of $\|\cdot\|$, $\|\cdot\|_{\pi}$, $\|\cdot\|_{F}$ and $\|\cdot\|_F^{\pi}$, and then it suffices to prove that 4) holds.
%
%

4) Referring to Lemma 1 in \cite{xin2019distributed}, one can obtain that (\ref{equ23}) and $\sigma=\|\mathcal{A}-\mathcal{A}_{\infty}\|_{\pi}<1$ hold. It suffices to prove (\ref{equ24}). Since
\begin{align*}
\mathcal{A}{\bf x}-\mathcal{A}_{\infty}{\bf x}=(\mathcal{A}-\mathcal{A}_{\infty})({\bf x}-\mathcal{A}_{\infty}{\bf x}),
\end{align*}
where ${\mathcal{A}\mathcal{A}_{\infty}}=\mathcal{A}_{\infty}\mathcal{A}=\mathcal{A}_{\infty}\mathcal{A}_{\infty}=\mathcal{A}_{\infty}$ is used, one has that
\begin{align*}
&\|\mathcal{A}{\bf x}-\mathcal{A}_{\infty}{\bf x}\|_{F}^{\pi}\\
&=\|{\rm diag}(\sqrt{\pi})(\mathcal{A}-\mathcal{A}_{\infty})({\bf x}-\mathcal{A}_{\infty}{\bf x})\|_{F}\\
&=\|{\rm diag}(\sqrt{\pi})(\mathcal{A}-\mathcal{A}_{\infty}){\rm diag}(\sqrt{\pi})^{-1}\\&~~~\times{\rm diag}(\sqrt{\pi})({\bf x}-\mathcal{A}_{\infty}{\bf x})\|_{F}\\
&\leq\|{\rm diag}(\sqrt{\pi})(\mathcal{A}-\mathcal{A}_{\infty}){\rm diag}(\sqrt{\pi})^{-1}\|\\&~~~\times\|{\rm diag}(\sqrt{\pi})({\bf x}-\mathcal{A}_{\infty}{\bf x})\|_{F}\\
&=\|\mathcal{A}-\mathcal{A}_{\infty}\|_{\pi}\|{\bf x}-\mathcal{A}_{\infty}{\bf x}\|_{F}^{\pi}\\
&=\sigma\|{\bf x}-\mathcal{A}_{\infty}{\bf x}\|_{F}^{\pi}.
\end{align*}
where the inequality is obtained based on 1) of this lemma.
The proof is completed. \hfill$\blacksquare$

\section{Proof of Theorem \ref{thm1}:}\label{appendix_C}
To prove Theorem \ref{thm1}, if we can prove the following two inequalities hold:
\begin{align}
&\|{\bf x}^{t+1}-\tilde{\bf x}^*\|_F^{\pi}\nonumber\\
&\leq\sqrt{\rho(M_{\alpha})}\|{\bf x}^t-\tilde{\bf x}^*\|_F^{\pi}+\frac{\sqrt{2}\alpha}{2\sqrt{n+m}}\sum_{i=1}^m\|{\bf v}_i^t-\overline{\bf v}_i^t\|_F,\label{equ-A1}\\
&\sum_{i=1}^m\|{\bf v}_i^{t+1}-\overline{\bf v}_i^{t+1}\|_F\nonumber\\
&\leq\sigma_{\max}\sum_{i=1}^m\|{\bf v}_i^t-\overline{\bf v}_i^t\|_F+\sqrt{m(n+m)}L\|{\bf x}^{t+1}-\tilde{\bf x}^*\|_F^{\pi}\nonumber\\
&~~~+\sqrt{m(n+m)}L\|{\bf x}^{t}-\tilde{\bf x}^*\|_F^{\pi},\label{equ-A2}
\end{align}
then it is obtained that
\begin{align}
\left[
\begin{array}{c}
\|{\bf x}^{t+1}-\tilde{\bf x}^*\|_F^{\pi}\\
\sum\limits_{i=1}^m\|{\bf v}_i^{t+1}-\overline{\bf v}_i^{t+1}\|_F
\end{array}
\right]
\leq H_{\alpha}
\left[
\begin{array}{c}
\|{\bf x}^{t}-\tilde{\bf x}^*\|_F^{\pi}\\
\sum\limits_{i=1}^m\|{\bf v}_i^{t}-\overline{\bf v}_i^{t}\|_F
\end{array}
\right].
\end{align}
By appropriately choosing $\alpha$ such that $\rho(M_{\alpha})<1$ and $\rho(H_{\alpha})<1$, it can be concluded that ${\bf x}^t$ converges to $\tilde{\bf x}^*$ linearly, and the proof of Theorem \ref{thm1} is thus completed.

In the following, we prove (\ref{equ-A1}) and (\ref{equ-A2}) are satisfied.

{\bf Proof of (\ref{equ-A1})}: In view of iteration (\ref{equ14}), one has
\begin{align}
&\|{\bf x}^{t+1}-\tilde{\bf x}^*\|_F^{\pi}\nonumber\\
&=\|\mathcal{P}_{\Omega}\left[\mathcal{A}{\bf x}^t-\Lambda V^t\right]-\mathcal{P}_{\Omega}\left[\tilde{\bf x}^*-\Lambda {\bf F}(\tilde{x}^*)\right]\|_F^{\pi}\nonumber\\
&\leq\|\mathcal{A}{\bf x}^t-\tilde{\bf x}^*-\Lambda V^t+\Lambda {\bf F}(\tilde{x}^*)\|_F^{\pi}\nonumber\\
&\leq\|\mathcal{A}{\bf x}^t-\tilde{\bf x}^*-\Lambda(\overline{V}^t-{\bf F}(\tilde{x}^*))\|_F^{\pi}+\|\Lambda({V}^t-\overline{V}^t)\|_F^{\pi},\label{equ-c4}
\end{align}
where
\begin{align}
{\bf F}(\tilde{\bf x}^*)&:={\rm diag}\left\{\frac{{\bf1}_{n_i}}{n_i}\sum_{j=1}^{n_i}(\nabla_if_{ij}(\tilde{x}^*))^{\top},i\in[m]\right\},\\
\overline{V}^t&:={\rm diag}\{\overline{\bf v}_i^t,i\in[m]\},
\end{align}
the equality is obtained by $\tilde{\bf x}^*=\mathcal{P}_{\Omega}\left[\tilde{\bf x}^*-\Lambda {\bf F}(\tilde{x}^*)\right]$ following (\ref{equ9e}) and the definition of $\mathcal{P}_{\Omega}[\cdot]$ in (\ref{e19}), and the first inequality is derived by
\begin{align*}
&\|\mathcal{P}_{\Omega}[A]-\mathcal{P}_{\Omega}[B]\|_F^{\pi}\\
&=\sqrt{\left\langle{\rm diag}(\pi)(\mathcal{P}_{\Omega}[A]-\mathcal{P}_{\Omega}[B]),\mathcal{P}_{\Omega}[A]-\mathcal{P}_{\Omega}[B]\right\rangle}\\
&=\sqrt{\sum_{h=1}^n\pi_h\|P_{\bm{\Omega}_h}[Col_h(A^{\top})]-P_{\bm{\Omega}_h}[Col_h(B^{\top})]\|^2}\\
&\leq\sqrt{\sum_{h=1}^n\pi_h\|Col_h(A^{\top})-Col_h(B^{\top})\|^2}\\
&=\|A-B\|_F^{\pi}, ~\forall A,B\in\mathbb{R}^{n\times q},
\end{align*}
relying on the definitions of $\mathcal{P}_{\Omega}[\cdot]$ in (\ref{e19}), the norm $\|\cdot\|_F^{\pi}$ in (\ref{equ19}) and the property of projection $P_{\Omega}[\cdot]$.

Decompose ${\bf x}^t$ as ${\bf x}^t=\mathcal{A}_{\infty}{\bf x}^t+{\bf x}^t_{\bot}$, where ${\bf x}^t_{\bot}:=(I_{n}-\mathcal{A}_{\infty}){\bf x}^t$ and $\mathcal{A}_{\infty}={\bf1}_n\pi^{\top}$. Then,
\begin{align}
\mathcal{A}{\bf x}^t&=\mathcal{A}\mathcal{A}_{\infty}{\bf x}^t+\mathcal{A}{\bf x}^t_{\bot}\nonumber\\
&=\mathcal{A}_{\infty}{\bf x}^t+\mathcal{A}{\bf x}^t_{\bot}.\label{equ-c7}
\end{align}
It can be verified that
\begin{align}
&\left\langle{\rm diag}(\pi)(\mathcal{A}_{\infty}{\bf x}^t-\tilde{\bf x}^*),{\bf x}^t_{\bot}\right\rangle\nonumber\\
&={\rm trace}[({\bf x}^t-\tilde{\bf x}^*)^{\top}\mathcal{A}_{\infty}^{\top}{\rm diag}(\pi)(I_{n}-\mathcal{A}_{\infty}){\bf x}^t]\nonumber\\
&=0,
\end{align}
where $\mathcal{A}_{\infty}\tilde{\bf x}^*=\tilde{\bf x}^*$ is used.
Similarly
\begin{align}
\left\langle{\rm diag}(\pi)\mathcal{A}_{\infty}{\bf x}^t,{\bf x}^t_{\bot}\right\rangle&=0,\\
\left\langle{\rm diag}(\pi)(\mathcal{A}_{\infty}{\bf x}^t-\tilde{\bf x}^*),\mathcal{A}{\bf x}^t_{\bot}\right\rangle&=0.\label{equ-c9}
\end{align}
Then, for the first term of the right-hand side of (\ref{equ-c4}), one has
\begin{align}
&\left(\|\mathcal{A}{\bf x}^t-\tilde{\bf x}^*-\Lambda(\overline{V}^t-{\bf F}(\tilde{\bf x}^*))\|_F^{\pi}\right)^2\nonumber\\
&=\left(\|\mathcal{A}_{\infty}{\bf x}^t-\tilde{\bf x}^*+\mathcal{A}{\bf x}^t_{\bot}-\Lambda(\overline{V}^t-{\bf F}(\tilde{\bf x}^*))\|_F^{\pi}\right)^2\nonumber\\
&=\left(\|\mathcal{A}_{\infty}{\bf x}^t-\tilde{\bf x}^*\|_F^{\pi}\right)^2+\left(\|\mathcal{A}{\bf x}^t_{\bot}\|_F^{\pi}\right)^2+\left(\|\Lambda(\overline{V}^t-{\bf F}(\tilde{\bf x}^*))\|_F^{\pi}\right)^2\nonumber\\
&~~~-2\left\langle{\rm diag}(\pi)(\mathcal{A}_{\infty}{\bf x}^t-\tilde{\bf x}^*),\Lambda(\overline{V}^t-{\bf F}(\tilde{\bf x}^*))\right\rangle\nonumber\\
&~~~-2\left\langle{\rm diag}(\pi)\mathcal{A}{\bf x}^t_{\bot},\Lambda(\overline{V}^t-{\bf F}(\tilde{\bf x}^*))\right\rangle,\label{equ-c10}
\end{align}
where (\ref{equ-c7}) is used to derive the first equality and (\ref{equ-c9}) is applied to obtain the second equality. To derive an upper bound on (\ref{equ-c10}), we analyze each term in (\ref{equ-c10}). Note that ${\bf x}^t_{\bot}=(I_{n}-\mathcal{A}_{\infty}){\bf x}^t$ and $\mathcal{A}\mathcal{A}_{\infty}=\mathcal{A}_{\infty}$, then based on (\ref{equ24}) in Lemma \ref{lemma1}, one has
\begin{align}\label{equ-c11}
\|\mathcal{A}{\bf x}^t_{\bot}\|_F^{\pi}
&=\|(\mathcal{A}-\mathcal{A}_{\infty}){\bf x}^t\|_F^{\pi}\nonumber\\
&\leq\sigma\|{\bf x}^t-\mathcal{A}_{\infty}{\bf x}^t\|_F^{\pi}\nonumber\\
&=\sigma\|{\bf x}^t_{\bot}\|_F^{\pi}.
\end{align}
In view of that $A_i\in\mathbb{R}^{n_i\times n_i}$ is a row and column stochastic matrix, it follows from (\ref{equ15}) that
\begin{align*}
{\bf1}_{n_i}^{\top}{\bf v}_i^{t+1}&={\bf1}_{n_i}^{\top}{\bf v}_i^t+\sum\limits_{j=1}^{n_i}\left[\nabla_if_{ij}(x_{ij}^{t+1},x_{-(ij)}^{t+1})\right]^{\top}\nonumber\\
&~~~-\sum\limits_{j=1}^{n_i}\left[\nabla_if_{ij}(x_{ij}^{t},x_{-(ij)}^{t})\right]^{\top},
\end{align*}
which indicates that for any $t\geq0$,
\begin{align*}
&{\bf1}_{n_i}^{\top}{\bf v}_i^{t+1}-\sum\limits_{j=1}^{n_i}\left[\nabla_if_{ij}(x_{ij}^{t+1},x_{-(ij)}^{t+1})\right]^{\top}\\
&={\bf1}_{n_i}^{\top}{\bf v}_i^t-\sum\limits_{j=1}^{n_i}\left[\nabla_if_{ij}(x_{ij}^{t},x_{-(ij)}^{t})\right]^{\top}.
\end{align*}
Therefore, one has
\begin{align}\label{equ36}
&{\bf1}_{n_i}^{\top}{\bf v}_i^{t}-\sum\limits_{j=1}^{n_i}\left[\nabla_if_{ij}(x_{ij}^{t},x_{-(ij)}^{t})\right]^{\top}={\bf0}_{1\times q_i},~\forall t\geq0,
\end{align}
since the initial value $v_{ij}^0=\nabla_if_{ij}(x_{ij}^0,x_{-(ij)}^0)$, $j\in[n_i]$. As a consequence, one can obtain that
\begin{align}
&\|\Lambda(\overline{V}^t-{\bf F}(\tilde{\bf x}^*))\|_F^{\pi}\nonumber\\
&\leq\frac{\sqrt{2}}{\sqrt{n+m}}\|\Lambda\|\|\overline{V}^t-{\bf F}(\tilde{\bf x}^*)\|_F\nonumber\\
&\leq\frac{\sqrt{2}\alpha}{2\sqrt{n+m}}\sqrt{\sum_{i=1}^m\left\|\frac{{\bf1}_{n_i}}{n_i}\sum_{j=1}^{n_i}[\nabla_if_{ij}(x_{ij}^t,x_{-(ij)}^t)-\nabla_if_{ij}(\tilde{x}^*)]^{\top}\right\|_F^2}\nonumber\\
&=\frac{\sqrt{2}\alpha}{2\sqrt{n+m}}\sqrt{\sum_{i=1}^m\frac{1}{n_i}\left\|\sum_{j=1}^{n_i}[\nabla_if_{ij}(x_{ij}^t,x_{-(ij)}^t)-\nabla_if_{ij}(\tilde{x}^*)]\right\|^2}\nonumber\\
&\leq\frac{\sqrt{2}\alpha}{2\sqrt{n+m}}\sqrt{\sum_{i=1}^m\sum_{j=1}^{n_i}\left\|\nabla_if_{ij}(x_{ij}^t,x_{-(ij)}^t)-\nabla_if_{ij}(\tilde{x}^*)\right\|^2}\nonumber\\
&\leq\frac{\sqrt{2}\alpha}{2\sqrt{n+m}}\sqrt{\sum_{i=1}^m\sum_{j=1}^{n_i}L_{ij}\left\|x_{(ij)}^t-\tilde{x}^*\right\|^2}\nonumber\\
&\leq\frac{\sqrt{2}L\alpha}{2\sqrt{n+m}}\|{\bf x}^t-\tilde{\bf x}^*\|_F\nonumber\\
&\leq\frac{\sqrt{2}L\alpha}{2}\|{\bf x}^t-\tilde{\bf x}^*\|_F^{\pi},\label{equ-c13}
\end{align}
where the first inequality is based on 1) and 3) in Lemma \ref{lemma1}, the second inequality is obtained by $\|\Lambda\|=\|{\rm diag}\{\alpha/(n_i+1),i\in[m]\}\|\leq\alpha/2$ and (\ref{equ36}), the fourth inequality holds under Assumption \ref{assumption2}, and the last inequality relies on 3) in Lemma \ref{lemma1}.

On the other hand, denote $x_{\pi}^t:=(\pi^{\top}{\bf x}^t)^{\top}\in\mathbb{R}^{q}$, then it can be derived that
\begin{align}
&-2\left\langle{\rm diag}(\pi)(\mathcal{A}_{\infty}{\bf x}^t-\tilde{\bf x}^*),\Lambda(\overline{V}^t-{\bf F}(\tilde{\bf x}^*))\right\rangle\nonumber\\
&=-2\left\langle{\rm diag}(\pi)(\mathcal{A}_{\infty}{\bf x}^t-\tilde{\bf x}^*),\Lambda(\overline{V}^t-{\bf F}(\mathcal{A}_{\infty}{\bf x}^t))\right\rangle\nonumber\\
&~~~-2\left\langle{\rm diag}(\pi)(\mathcal{A}_{\infty}{\bf x}^t-\tilde{\bf x}^*),\Lambda({\bf F}(\mathcal{A}_{\infty}{\bf x}^t)-{\bf F}(\tilde{\bf x}^*))\right\rangle\nonumber\\
&\leq2\|\mathcal{A}_{\infty}{\bf x}^t-\tilde{\bf x}^*\|_F^{\pi}\cdot\|\Lambda(\overline{V}^t-{\bf F}(\mathcal{A}_{\infty}{\bf x}^t))\|_F^{\pi}\nonumber\\
&~~~-2\left\langle \pi(x_{\pi}^t-\tilde{x}^*)^{\top},\Lambda({\bf F}(\mathcal{A}_{\infty}{\bf x}^t)-{\bf F}(\tilde{\bf x}^*))\right\rangle\nonumber\\
&=2\|\mathcal{A}_{\infty}{\bf x}^t-\tilde{\bf x}^*\|_F^{\pi}\cdot\|\Lambda(\overline{V}^t-{\bf F}(\mathcal{A}_{\infty}{\bf x}^t))\|_F^{\pi}\nonumber\\
&~~~-\frac{2\alpha}{n+m}(F(x_{\pi}^t)-F(\tilde{x}^*))^{\top}(x_{\pi}^t-\tilde{x}^*)\nonumber\\
&\leq\sqrt{2}L\alpha\|\mathcal{A}_{\infty}{\bf x}^t-\tilde{\bf x}^*\|_F^{\pi}\cdot\|{\bf x}^t_{\bot}\|_F^{\pi}-\frac{2\mu\alpha}{n+m}\|x_{\pi}^t-\tilde{x}^*\|^2\nonumber\\
&\leq\sqrt{2}L\alpha\|\mathcal{A}_{\infty}{\bf x}^t-\tilde{\bf x}^*\|_F^{\pi}\cdot\|{\bf x}_{\bot}^t\|_F^{\pi}\nonumber\\
&~~~-\frac{2\mu\alpha}{n}\left(\|\mathcal{A}_{\infty}{\bf x}^t-\tilde{\bf x}^*\|_F^{\pi}\right)^2,
\end{align}
where
\begin{align*}
{\bf F}(\mathcal{A}_{\infty}{\bf x}^t):={\rm diag}\left\{\frac{{\bf1}_{n_i}}{n_i}\sum_{j=1}^{n_i}(\nabla_if_{ij}(x_{\pi}^t))^{\top},i\in[m]\right\},
\end{align*}
$F(x_{\pi}^t)$ and $F(\tilde{x}^*)$ are defined in (\ref{equ5}) by letting $y=x_{\pi}^t$ and $y=\tilde{x}^*$, respectively, the second equality is due to Lemma \ref{lem1}, and the second inequality is obtained based on Assumption \ref{assumption3} and the following inequality
\begin{align}
\|\Lambda(\overline{V}^t-{\bf F}(\mathcal{A}_{\infty}{\bf x}^t))\|_F^{\pi}&\leq\frac{\sqrt{2}L\alpha}{2}\|{\bf x}^t-\mathcal{A}_{\infty}{\bf x}^t\|_F^{\pi}\nonumber\\
&=\frac{\sqrt{2}L\alpha}{2}\|{\bf x}_{\bot}^t\|_F^{\pi}
\end{align}
based on (\ref{equ-c13}).

Also, it can be obtained from (\ref{equ-c13}) that
\begin{align}
&-2\left\langle{\rm diag}(\pi)\mathcal{A}{\bf x}^t_{\bot},\Lambda(\overline{V}^t-{\bf F}(\tilde{\bf x}^*))\right\rangle\nonumber\\
&\leq2\|\mathcal{A}{\bf x}^t_{\bot}\|_F^{\pi}\cdot\|\Lambda(\overline{V}^t-{\bf F}(\tilde{\bf x}^*))\|_F^{\pi}\nonumber\\
&\leq\sqrt{2}L\alpha\|\mathcal{A}{\bf x}^t_{\bot}\|_F^{\pi}\cdot\|{\bf x}^t-\tilde{\bf x}^*\|_F^{\pi}\nonumber\\
&\leq\sqrt{2}L\alpha\sigma\|{\bf x}^t_{\bot}\|_F^{\pi}\cdot(\|\mathcal{A}_{\infty}{\bf x}^t-\tilde{\bf x}^*\|_F^{\pi}+\|{\bf x}^t_{\bot}\|_F^{\pi}).\label{equ-c16}
\end{align}

Therefore, substituting (\ref{equ-c11})--(\ref{equ-c16}) into (\ref{equ-c10}) yields that
\begin{align}
&\left(\|\mathcal{A}{\bf x}^t-\tilde{\bf x}^*-\Lambda(\overline{V}^t-{\bf F}(\tilde{\bf x}^*))\|_F^{\pi}\right)^2\nonumber\\
&\leq\left(\|\mathcal{A}_{\infty}{\bf x}^t-\tilde{\bf x}^*\|_F^{\pi},\|{\bf x}^t_{\bot}\|_F^{\pi}\right)
  M_{\alpha}
  \left(\begin{array}{c}\|\mathcal{A}_{\infty}{\bf x}^t-\tilde{\bf x}^*\|_F^{\pi}\\
  \|{\bf x}^t_{\bot}\|_F^{\pi}
  \end{array}
  \right)\nonumber\\
&\leq\rho(M_{\alpha})\left((\|\mathcal{A}_{\infty}{\bf x}^t-\tilde{\bf x}^*\|_F^{\pi})^2+(\|{\bf x}^t_{\bot}\|_F^{\pi})^2\right)\nonumber\\
&=\rho(M_{\alpha})\left(\|{\bf x}^t-\tilde{\bf x}^*\|_F^{\pi}\right)^2,\label{equ-c17}
\end{align}
where $(\|{\bf x}^t-\tilde{\bf x}^*\|_F^{\pi})^2=(\|\mathcal{A}_{\infty}{\bf x}^t-\tilde{\bf x}^*\|_F^{\pi})^2+(\|{\bf x}_{\bot}^t\|_F^{\pi})^2$ is applied to derive the first inequality.
Then, it can be obtained from (\ref{equ-c17}) that
\begin{align}
&\|\mathcal{A}{\bf x}^t-\tilde{\bf x}^*-\Lambda(\overline{V}^t-{\bf F}(\tilde{\bf x}^*))\|_F^{\pi}\nonumber\\
&\leq\sqrt{\rho(M_{\alpha})}\|{\bf x}^t-\tilde{\bf x}^*\|_F^{\pi}.\label{equ-c18}
\end{align}

For the second term in the right-hand side of (\ref{equ-c4}), it holds
\begin{align}
\|\Lambda({V}^t-\overline{V}^t)\|_F^{\pi}
&\leq\frac{\sqrt{2}}{\sqrt{n+m}}\|\Lambda\|\|{V}^t-\overline{V}^t\|_F\nonumber\\
&\leq\frac{\sqrt{2}\alpha}{2\sqrt{n+m}}\sqrt{\sum_{i=1}^m\|{\bf v}_i^t-\overline{\bf v}_i^t\|_F^2}\nonumber\\
&\leq\frac{\sqrt{2}\alpha}{2\sqrt{n+m}}\sum_{i=1}^m\|{\bf v}_i^t-\overline{\bf v}_i^t\|_F,\label{equ-c19}
\end{align}
where $\|\Lambda\|\leq\alpha/2$ is used to obtain the second inequality.

By (\ref{equ-c18}) and (\ref{equ-c19}), it can be derived from (\ref{equ-c4}) that
(\ref{equ-A1}) is satisfied.

{\bf Proof of (\ref{equ-A2}):}
For $\sum_{i=1}^m\|{\bf v}_i^t-\overline{\bf v}_i^t\|_{F}$, considering iteration (\ref{equ15}), the iteration of $\overline{\bf v}_i^t$ is obtained as
\begin{align}
\overline{\bf v}_i^{t+1}=\overline{\bf v}_i^t+\frac{{\bf1}_{n_i}{\bf1}_{n_i}^{\top}}{n_i}G_i({\bf x}_i^{t+1})-\frac{{\bf1}_{n_i}{\bf1}_{n_i}^{\top}}{n_i}G_i({\bf x}_i^{t}),
\end{align}
where ${\bf1}_{n_i}^{\top}A_i={\bf1}_{n_i}^{\top}$ is used. Then we obtain
\begin{align}
&\|{\bf v}_i^{t+1}-\overline{\bf v}_i^{t+1}\|_{F}\nonumber\\
&=\|A_i{\bf v}_i^{t}+G_i({\bf x}^{t+1})-G_i({\bf x}^{t})-\overline{\bf v}_i^{t}-\frac{{\bf1}_{n_i}{\bf1}_{n_i}^{\top}}{n_i}G_i({\bf x}_i^{t+1})\nonumber\\
&~~~~~+\frac{{\bf1}_{n_i}{\bf1}_{n_i}^{\top}}{n_i}G_i({\bf x}_i^{t})\|_{F}\nonumber\\
&\leq\|A_i{\bf v}_i^{t}-\overline{\bf v}_i^{t}\|_{F}+\|I_{n_i}-{\bf1}_{n_i}{\bf1}_{n_i}^{\top}/n_i\|\|G_i({\bf x}^{t+1})-G_i({\bf x}^{t})\|_{F}\nonumber\\
&=\|(A_i-{\bf1}_{n_i}{\bf1}_{n_i}^{\top}/{n_i})({\bf v}_i^{t}-\overline{\bf v}_i^{t})\|_{F}+\|G_i({\bf x}^{t+1})-G_i({\bf x}^{t})\|_{F}\nonumber\\
&\leq\|A_i-{\bf1}_{n_i}{\bf1}_{n_i}^{\top}/{n_i}\|\|{\bf v}_i^{t}-\overline{\bf v}_i^{t}\|_{F}\nonumber\\
&~~~+\sqrt{\sum\limits_{j=1}^{n_i}\left\|\nabla_if_{ij}(x_{(ij)}^{t+1})-\nabla_if_{ij}(x_{(ij)}^{t})\right\|^2}\nonumber\\
&\leq\sigma_i\|{\bf v}_i^{t}-\overline{\bf v}_i^{t}\|_{F}+\sqrt{\sum\limits_{j=1}^{n_i}L_{ij}^2\left\|x_{(ij)}^{t+1}-x_{(ij)}^{t}\right\|^2}\nonumber\\
&\leq\sigma_i\|{\bf v}_i^{t}-\overline{\bf v}_i^{t}\|_{F}+L\|{\bf x}_{i}^{t+1}-{\bf x}_{i}^{t}\|_{F},
\end{align}
where the second equality is obtained based on $A_i{\bf1}_{n_i}={\bf1}_{n_i}$ and $\|I_{n_i}-{\bf1}_{n_i}{\bf1}_{n_i}^{\top}/{n_i}\|=1$, the second inequality is based on 1) in Lemma \ref{lemma1} and (\ref{e18}), the third inequality is derived according to (\ref{e3}) and Assumption \ref{assumption2}, and the last inequality hinges on the structure of ${\bf x}^t_i$ in (\ref{e14}). Therefore, one has that
\begin{align}
&\sum\limits_{i=1}^m\|{\bf v}_i^{t+1}-\overline{\bf v}_i^{t+1}\|_{F}\nonumber\\
&\leq\sum\limits_{i=1}^m\sigma_i\|{\bf v}_i^{t}-\overline{\bf v}_i^{t}\|_{F}+L\sum\limits_{i=1}^m\|{\bf x}_{i}^{t+1}-{\bf x}_{i}^{t}\|_{F}\nonumber\\
&\leq\sigma_{\rm max}\sum\limits_{i=1}^m\|{\bf v}_i^{t}-\overline{\bf v}_i^{t}\|_{F}+L\sqrt{m}\|{\bf x}^{t+1}-{\bf x}^{t}\|_{F},\nonumber\\
&\leq\sigma_{\rm max}\sum\limits_{i=1}^m\|{\bf v}_i^{t}-\overline{\bf v}_i^{t}\|_{F}\nonumber\\
&~~~+L\sqrt{m}\sqrt{n+m}(\|{\bf x}^{t+1}-\tilde{\bf x}^{*}\|_{F}^{\pi}+\|{\bf x}^{t}-\tilde{\bf x}^{*}\|_{F}^{\pi}),\label{equ53}
\end{align}
where $\sigma_{\max}=\max_{i\in[m]}\sigma_i$. Therefore, (\ref{equ-A2}) is satisfied. \hfill$\blacksquare$

\section{Proof of Proposition \ref{pro1}:}\label{appendix_D}
Note that $M_{\alpha}$ is symmetric and each entry of $M_{\alpha}$ is positive by letting $(M_{\alpha})_{1,1}>0$, ensured by $\alpha<\frac{n}{2\mu}$. Then, $\rho(M_{\alpha})<1$ holds if $0\prec M_{\alpha}\prec I_2$. By the diagonal dominance, $M_{\alpha}\succ0$ if $(M_{\alpha})_{1,1}-(M_{\alpha})_{1,2}>0$ and $(M_{\alpha})_{2,2}-(M_{\alpha})_{2,1}>0$, that is,
\begin{align}
1-\frac{2\mu\alpha}{n}+\frac{L^2\alpha^2}{2}-\frac{\sqrt{2}(1+\sigma)L\alpha}{2}&>0,\label{equ-D1}\\
\sigma^2+\sqrt{2}\sigma L\alpha+\frac{L^2\alpha^2}{2}-\frac{\sqrt{2}(1+\sigma)L\alpha}{2}&>0,\label{equ-D2}
\end{align}
which can be ensured by
\begin{align*}
\frac{2\mu\alpha}{n}&<\frac{1}{2},\\
\frac{\sqrt{2}(1+\sigma)L\alpha}{2}&<\frac{1}{2},\\
\sigma^2+\sqrt{2}\sigma L\alpha-\frac{\sqrt{2}(1+\sigma)L\alpha}{2}&>0,
\end{align*}
i.e.,
\begin{align}\label{equ-D3}
\alpha<\min\left\{\frac{n}{4\mu},\frac{1}{\sqrt{2}(1+\sigma)L},\frac{\sqrt{2}\sigma^2}{(1-\sigma)L}\right\}.
\end{align}
Furthermore, by the Sylvester's criterion, $M_{\alpha}\prec I_2$ can be ensured if $(I_2-M_{\alpha})_{1,1}>0$ and $\det(I_2-M_{\alpha})>0$, that is,
\begin{align}
&\frac{2\mu\alpha}{n}-\frac{L^2\alpha^2}{2}>0,\label{equ-D4}\\
&\left(\frac{2\mu\alpha}{n}-\frac{L^2\alpha^2}{2}\right)\left(1-\sigma^2-\sqrt{2}\sigma L\alpha-\frac{L^2\alpha^2}{2}\right)\nonumber\\
&-\left(\frac{\sqrt{2}(1+\sigma)L\alpha}{2}\right)^2>0.\label{equ-D5}
\end{align}
It can be easily seen that (\ref{equ-D4}) is equivalent to
\begin{align}\label{equ-D6}
\alpha<\frac{4\mu}{nL^2}.
\end{align}
After a simple computation, (\ref{equ-D5}) is
\begin{align*}
&\frac{2\mu(1-\sigma^2)}{n}\alpha-\frac{2\sqrt{2}\mu\sigma L}{n}\alpha^2\\
&-(1+\sigma)L^2\alpha^2-\frac{\mu L^2}{n}\alpha^3+\frac{\sqrt{2}\sigma L^3}{2}\alpha^3+\frac{L^4}{4}\alpha^4>0,
\end{align*}
which is implied by
\begin{align*}
\frac{2\mu(1-\sigma^2)}{3n}\alpha-\frac{2\sqrt{2}\mu\sigma L}{n}\alpha^2&>0,\\
\frac{2\mu(1-\sigma^2)}{3n}\alpha-(1+\sigma)L^2\alpha^2&>0,\\
\frac{2\mu(1-\sigma^2)}{3n}\alpha-\frac{\mu L^2}{n}\alpha^3&>0,
\end{align*}
that is,
\begin{align}\label{equ-D7}
\alpha<\min\left\{\frac{1-\sigma^2}{3\sqrt{2}\sigma L},\frac{2\mu(1-\sigma)}{3nL^2},\frac{\sqrt{2(1-\sigma^2)}}{\sqrt{3}L}\right\}.
\end{align}
Combining (\ref{equ-D3}), (\ref{equ-D6}) and (\ref{equ-D7}) yields that if
\begin{align}
&\alpha<\min\left\{\frac{n}{4\mu},\frac{1}{\sqrt{2}(1+\sigma)L},\frac{\sqrt{2}\sigma^2}{(1-\sigma)L},\right.\nonumber\\
&\left.~~~~~~~~~~~~\frac{1-\sigma^2}{3\sqrt{2}\sigma L},\frac{2\mu(1-\sigma)}{3nL^2},\frac{\sqrt{2(1-\sigma^2)}}{\sqrt{3}L}\right\},\label{equ-D8}
\end{align}
then $\rho(M_{\alpha})<1$.

On the other hand, based on Corollary 8.1.29 in \cite{horn2012}, if there exists a vector $v=col(v_1,v_2)$ with $v_1>0,v_2>0$ such that $H_{\alpha}v<v$, then $\rho(H_{\alpha})<1$. Note that $H_{\alpha}v<v$ is equivalent to
\begin{align}
&\sqrt{\rho(M_{\alpha})}v_1+\frac{\sqrt{2}\alpha}{2\sqrt{n+m}}v_2<v_1,\label{equ-D9}\\
&\sqrt{m(n+m)}(1+\sqrt{\rho(M_{\alpha})})Lv_1+(\sigma_{\max}+\frac{\sqrt{2m}L\alpha}{2})v_2<v_2.\label{equ-D10}
\end{align}
Let $v_2=1$, then it can be obtained from (\ref{equ-D9}) that
\begin{align}
\alpha<\sqrt{2(n+m)}(1-\sqrt{\rho(M_{\alpha})})v_1.\label{equ-D11}
\end{align}
By taking $v_1=\frac{1-\sigma_{\max}}{2\sqrt{m(n+m)}(1+\sqrt{\rho(M_{\alpha})})L}$, then one has from (\ref{equ-D10}) and (\ref{equ-D11}) that
\begin{align}\label{equ-D12}
\alpha<\left\{\frac{1-\sigma_{\max}}{\sqrt{2m}L},\frac{(1-\sigma_{\max})(1-\sqrt{\rho(M_{\alpha})})}{\sqrt{2m}(1+\sqrt{\rho(M_{\alpha})})L}\right\}.
\end{align}
Combining (\ref{equ-D8}) and (\ref{equ-D12}) yields (\ref{e32}). \hfill$\blacksquare$

\bibliographystyle{elsarticle-num}


\begin{thebibliography}{10}
\expandafter\ifx\csname url\endcsname\relax
  \def\url#1{\texttt{#1}}\fi
\expandafter\ifx\csname urlprefix\endcsname\relax\def\urlprefix{URL }\fi
\expandafter\ifx\csname href\endcsname\relax
  \def\href#1#2{#2} \def\path#1{#1}\fi

\bibitem{li2015demand}
Li, N., Chen, L., \& Dahleh, M.~A. (2015). Demand response using linear supply function
  bidding. {\em IEEE Transactions on Smart Grid}, {\em 6}(4), 1827--1838.

\bibitem{saad2012game}
Saad, W., Han, Z., Poor, H. V., \& Basar, T. (2012). Game-theoretic methods for the smart
  grid: an overview of microgrid systems, demand-side management, and smart
  grid communications. {\em IEEE Signal Processing Magazine}, {\em 29}(5), 86--105.

\bibitem{hollander2006applicability}
Hollander, Y., \& Prashker, J. N. (2006). The applicability of non-cooperative game theory
  in transport analysis. {\em Transportation}, {\em 33}(5), 481--496.

\bibitem{facchinei2010generalized}
Facchinei, F., \& Kanzow, C. (2010). {Generalized Nash equilibrium problems}. {\em Annals of
  Operations Research}, {\em 175}(1), 177--211.

\bibitem{yu2017distributed}
Yu, C.~K., Van Der~Schaar, M., \& Sayed, A. H. (2017). {Distributed learning for stochastic
  generalized Nash equilibrium problems}. {\em IEEE Transactions on Signal
  Processing}, {\em 65}(15), 3893--3908.

\bibitem{shamma2005dynamic}
Shamma, J.~S., \& Arslan, G. (2005). {Dynamic fictitious play, dynamic gradient play, and
  distributed convergence to Nash equilibria}. {\em IEEE Transactions on Automatic
  Control}, {\em 50}(3), 312--327.

\bibitem{stankovic2012distributed}
Stankovic, M. S., Johansson, K. H., \& Stipanovic, D. M. (2012). {Distributed seeking of
  Nash equilibria with applications to mobile sensor networks}. {\em IEEE
  Transactions on Automatic Control}, {\em 57}(4), 904--919.

\bibitem{frihauf2011nash}
Frihauf, P., Krstic, M., \& Basar, T. (2011). Nash equilibrium seeking in noncooperative
  games. {\em IEEE Transactions on Automatic Control}, {\em 57}(5), 1192--1207.

\bibitem{de2019distributed}
De~Persis, C., \& Grammatico, S. (2019). {Distributed averaging integral Nash equilibrium
  seeking on networks}. {\em Automatica}, {\em 110}, 108548.

\bibitem{gadjov2019passivity}
Gadjov, D., \& Pavel, L. (2019). {A passivity-based approach to Nash equilibrium seeking
  over networks}. {\em IEEE Transactions on Automatic Control}, {\em 64}(3), 1077--1092.

\bibitem{koshal2016distributed}
Koshal, J., Nedi{\'c}, A., \& Shanbhag, U.~V. (2016). Distributed algorithms for aggregative
  games on graphs. {\em Operations Research}, {\em 64}(3), 680--704.

\bibitem{salehisadaghiani2016distributed}
Salehisadaghiani, F., \& Pavel, L. (2016). {Distributed Nash equilibrium seeking: a
  gossip-based algorithm}. {\em Automatica}, {\em 72}, 209--216.

\bibitem{salehisadaghiani2019distributed}
Salehisadaghiani, F., Shi, W., \& Pavel, L. (2019). {Distributed Nash equilibrium seeking
  under partial-decision information via the alternating direction method of
  multipliers}. {\em Automatica}, {\em 103}, 27--35.

\bibitem{tatarenko2019geometric}
Tatarenko,T., \& Nedich, A. (2019). Geometric convergence of distributed gradient play in
  games with unconstrained action sets. arXiv preprint arXiv:1907.07144.

\bibitem{tatarenko2018geometric}
Tatarenko, T., Shi, W., \& Nedich, A. (2021). Geometric convergence of gradient play
  algorithms for distributed Nash equilibrium seeking. {\em IEEE Transactions on Automatic Control}, {\em 66}(11), 5342--5353.

\bibitem{bianchi2021fully}
Bianchi, M., \& Grammatico, S. (2021). Fully distributed Nash equilibrium seeking over time-varying communication networks with linear convergence rate. {\em IEEE Control Systems Letters}, {\em 5}(2), 499-504.


\bibitem{nedic2009distributed}
Nedich, A., \& Ozdaglar, A. (2009). Distributed subgradient methods for multi-agent
  optimization. {\em IEEE Transactions on Automatic Control}, {\em 54}(1), 48--61.

\bibitem{zeng2017distributed}
Zeng, X., Yi, P., \& Hong, Y. (2017). Distributed continuous-time algorithm for constrained
  convex optimizations via nonsmooth analysis approach. {\em IEEE Transactions on
  Automatic Control}, {\em 62}(10), 5227--5233.

\bibitem{qu2018harnessing}
Qu, G., \& Li, N. (2018). Harnessing smoothness to accelerate distributed optimization. {\em IEEE Transactions on Control of Network Systems}, {\em 5}(3), 1245--1260.

\bibitem{li2020distri}
Li, X., Feng, G., \& Xie, L. (2021). Distributed proximal algorithms for multi-agent
  optimization with coupled inequality constraints. {\em IEEE Transactions on
  Automatic Control}, {\em 66}(3), 1223-1230.

\bibitem{li2019distributed_ijrnc}
Li, X., Xie, L., \& Hong, Y. (2019). Distributed continuous-time algorithm for a general
  nonsmooth monotropic optimization problem. {\em International Journal of Robust
  and Nonlinear Control}, {\em 29}(10), 3252--3266.

\bibitem{li2019distributed}
Li, X., Xie, L., \& Hong, Y. (2020). Distributed continuous-time nonsmooth convex
  optimization with coupled inequality constraints. {\em IEEE Transactions on
  Control of Network Systems}, {\em 7}(1), 74--84.

\bibitem{bullo2009distributed}
Bullo, F., Cortes, J., \& Martinez, S. (2009). {\em Distributed control of robotic networks: a
  mathematical approach to motion coordination algorithms}. Princeton University
  Press.

\bibitem{peng2009coexistence}
Peng, T.~J.~A., \& Bourne, M. (2009). The coexistence of competition and cooperation
  between networks: implications from two taiwanese healthcare networks,
  {\em British Journal of Management}, {\em 20}(3), 377--400.

\bibitem{shehory1998methods}
Shehory, O., \& Kraus, S. (1998). Methods for task allocation via agent coalition
  formation. {\em Artificial Intelligence}, {\em 101}(1-2), 165--200.

\bibitem{ye2018nash}
Ye, M., Hu, G., \& Lewis, F. L. (2018). {Nash equilibrium seeking for N-coalition
  noncooperative games}. {\em Automatica}, {\em 95}, 266--272.

\bibitem{ye2019unified}
Ye, M., Hu, G., Lewis, F. L., \& Xie, L. (2019). A unified strategy for solution seeking in
  graphical N-coalition noncooperative games. {\em IEEE Transactions on
  Automatic Control}, {\em 64}(11), 4645--4652.

\bibitem{ye2020extremum}
Ye, M., Hu, G., \& Xu, S. (2020). {An extremum seeking-based approach for Nash equilibrium
  seeking in N-cluster noncooperative games}. {\em Automatica}, {\em 114}, 108815.

\bibitem{zeng2019generalized}
Zeng, X., Chen, J., Liang, S., \& Hong, Y. (2019). Generalized nash equilibrium seeking
  strategy for distributed nonsmooth multi-cluster game. {\em Automatica}, {\em 103},
  20--26.

%
%


\bibitem{li2020}
Li, X., Xie, L., \& Hong, Y. (2021). Distributed aggregative optimization over multi-agent
  networks. {\em IEEE Transactions on Automatic Control}, DOI: 10.1109/TAC.2021.3095456.


\bibitem{xin2019distributed}
Xin, R., Sahu, A. K., Khan, U. A., \& Kar, S. (2019). Distributed stochastic optimization
  with gradient tracking over strongly-connected networks.  In \emph{2019 IEEE 58th Conference on Decision and Control (pp. 8353-8358).}


\bibitem{horn2012}
Horn, R. A., \& Johnson, C. R. (2012). {\em Matrix analysis}. Cambridge university press.


\end{thebibliography}

\end{document}